\documentclass[10pt]{amsart}
\usepackage[cp1251]{inputenc}
\usepackage[english]{babel}
\usepackage{amsmath}
\usepackage{amssymb}
\usepackage{amsfonts}

\begin{document}
\renewcommand{\refname}{References}

\thispagestyle{empty}

\title[To Numerical Modeling With Strong Orders 1.0, 1.5, and 2.0]
{To Numerical Modeling With Strong Orders 1.0, 1.5, and 2.0 of
Convergence for
Multidimensional Dynamical Systems With Random 
Disturbances}
\author[D.F. Kuznetsov]{Dmitriy F. Kuznetsov}
\address{Dmitriy Feliksovich Kuznetsov
\newline\hphantom{iii} Peter the Great Saint-Petersburg Polytechnic University,
\newline\hphantom{iii} Polytechnicheskaya ul., 29,
\newline\hphantom{iii} 195251, Saint-Petersburg, Russia}%
\email{sde\_kuznetsov@inbox.ru}

\maketitle {\small
\begin{quote}
\noindent{\sc Abstract.} 
The article is devoted to explicit one-step numerical methods
with strong orders 1.0, 1.5, and 2.0 of convergence for
Ito stochastic differential equations with
multidimensional and non-commutative noise.
For numerical modeling of iterated Ito
stochastic integrals with multiplicities 1 to 4 we use the method
of multiple Fourier--Legendre series converging in the sense
of norm 
in Hilbert space $L_2([t, T]^k),$ $k=1,2,3,4.$
The article is addressed to engineers who use
numerical mo\-deling in stochastic control
and for solving the nonlinear filtering problem.
\medskip
\end{quote}
}

\vspace{5mm}

%\linespread{1.6}

\setlength{\baselineskip}{2.0em}

\tableofcontents

\setlength{\baselineskip}{1.2em}

%\linespread{1.0}

\section{Introduction}

\vspace{5mm}

The Ito stochastic differential equations (SDEs) are known to be 
adequate mathematical models
of the dynamical systems of various physical nature subjected to random 
perturbations \cite{KlPl2}-\cite{Shir1}. On the
assumption of strong convergence criterion \cite{KlPl2}, 
the need for numerical 
integration of Ito SDEs
arises at solving the different mathematical problems.
Among them we mention the following problems:
stochastic optimal 
control (also with incomplete
data) \cite{KlPl2}, \cite{Lip}, 
signal filtering 
in random noise in
various formulations \cite{KlPl2}, \cite{Lip}, 
estimating the parameters of stochastic systems \cite{KlPl2}, \cite{KPS}.
It is common knowledge
that one of the promising
approaches to the numerical integration of Ito SDEs is the approach 
based on 
the stochastic analogues
of the Taylor formula, the so-called Taylor--Ito and Taylor--Stratonovich 
expansions \cite{KlPl2}, \cite{KPS}, \cite{Mi2}-\cite{kuz33}. This
approach makes use of finite discretization of the time variable 
and implies numerical modeling of
the solution of Ito SDE at the discrete time instants using the 
stochastic analogues of the Taylor
formula obtained by iterative application of the Ito formula.

Let $(\Omega,$ ${\rm F},$ ${\sf P})$ be a complete probability space, let 
$\{{\rm F}_t, t\in[0,T]\}$ be a nondecreasing right-continous 
family of $\sigma$-algebras of ${\rm F},$
and let ${\bf f}_t$ be a standard $m$-dimensional Wiener 
stochastic process, which is
${\rm F}_t$-measurable for any $t\in[0, T].$ We assume that the components
${\bf f}_{t}^{(i)}$ $(i=1,\ldots,m)$ of this process are 
independent. Consider
an Ito SDE in the integral form  

\begin{equation}
\label{1.5.2}
{\bf x}_t={\bf x}_0+\int\limits_0^t {\bf a}({\bf x}_{\tau},\tau)d\tau+
\int\limits_0^t B({\bf x}_{\tau},\tau)d{\bf f}_{\tau},\ \ \
{\bf x}_0={\bf x}(0,\omega),\ \ \ \omega\in\Omega.
\end{equation}

\vspace{3mm}
\noindent
Here ${\bf x}_t$ is some $n$-dimensional stochastic process 
satisfying to the Ito SDE (\ref{1.5.2}). 
The nonrandom functions ${\bf a}: \mathbb{R}^n\times[0, T]\to\mathbb{R}^n$,
$B: \mathbb{R}^n\times[0, T]\to\mathbb{R}^{n\times m}$
guarantee the existence and uniqueness up to stochastic 
equivalence of a solution
of the Ito SDE (\ref{1.5.2}) \cite{1}. The second integral 
on the right-hand side of (\ref{1.5.2}) is 
interpreted as an Ito stochastic integral.
Let ${\bf x}_0$ be an $n$-dimensional random variable, which is 
${\rm F}_0$-measurable and 
${\sf M}\{\left|{\bf x}_0\right|^2\}<\infty$ 
(${\sf M}$ denotes a mathematical expectation).
We assume that
${\bf x}_0$ and ${\bf f}_t-{\bf f}_0$ are independent when $t>0.$

The most important feature
of stochastic analogues of the Taylor formula 
\cite{KlPl2}, \cite{KPS}, \cite{Mi2}-\cite{kuz33} for solutions 
of the Ito SDE (\ref{1.5.2}) consists in the presence
of iterated Ito and Stratonovich stochastic integrals.
These stochastic integrals 
are complicated functionals from the components of the 
multidimensional Wiener process. In one of
the most general forms of notation of the present paper, the 
aforementioned iterated Ito and Stratonovich
stochastic integrals are given, respectively, by

\vspace{-1mm}
\begin{equation}
\label{ito}
J[\psi^{(k)}]_{T,t}=\int\limits_t^T\psi_k(t_k) \ldots \int\limits_t^{t_{2}}
\psi_1(t_1) d{\bf w}_{t_1}^{(i_1)}\ldots
d{\bf w}_{t_k}^{(i_k)},
\end{equation}

\begin{equation}
\label{str}
J^{*}[\psi^{(k)}]_{T,t}=
\int\limits_t^{*T}\psi_k(t_k) \ldots \int\limits_t^{*t_{2}}
\psi_1(t_1) d{\bf w}_{t_1}^{(i_1)}\ldots
d{\bf w}_{t_k}^{(i_k)},
\end{equation}

\vspace{3mm}
\noindent
where every $\psi_l(\tau)$ $(l=1,\ldots,k)$ is a 
nonrandom function 
on $[t,T],$ ${\bf w}_{\tau}^{(i)}={\bf f}_{\tau}^{(i)}$
for $i=1,\ldots,m$ and
${\bf w}_{\tau}^{(0)}=\tau,$\

\vspace{-3mm}
$$
\int\limits\ \hbox{and}\ \int\limits^{*}
$$ 

\vspace{3mm}
\noindent
denote Ito and 
Stratonovich stochastic integrals,
respectively,\ $i_1,\ldots,i_k = 0, 1,\ldots,m$
(in this paper, 
we use the definition of the Stratonovich stochastic integral from \cite{KlPl2}).

Consequently, the systems of stochastic integrals like 
(\ref{ito}), (\ref{str})
play an important part in
solving the problem of numerical integration of the Ito SDEs (\ref{1.5.2}). 
In terms of the mean-square
convergence criterion, the problem of efficient joint numerical 
modeling of the totalities of stochastic
integrals of the kind (\ref{ito}), (\ref{str}) (the case of a multidimensional 
Wiener process) is not only important,
but also sufficiently complex in both the theoretical and computational 
terms.
We note that the aforementioned problem does not arise 
at using the Euler method for the
Ito SDEs (\ref{1.5.2}) \cite{KlPl2}, \cite{Mi2}. 
However, despite its simplicity, the 
Euler method under the standard conditions
\cite{KlPl2}, \cite{Mi2} for coefficients of the Ito SDE (\ref{1.5.2}) 
has the 
mean-square convergence order 
0.5 \cite{KlPl2}, \cite{Mi2}, and its accuracy is insufficient to solve a 
number of practical problems. This fact
motivates one to construct numerical methods for the Ito SDEs 
(\ref{1.5.2}) having higher orders of strong
convergence.

It may seem at the first glance that the stochastic integrals 
from the families (\ref{ito}), (\ref{str})
can
be approximated by the multiple integral sums. However, 
this leads to partitioning of the interval
of integration $[t, T]$ of the iterated stochastic integrals.
The mentioned interval 
is already a small value because it
represents a step of integration in the numerical methods 
for Ito SDEs. As the numerical
experiments show \cite{7}, the above partitioning gives rise to an unacceptably 
high computing costs.

A number of publications are devoted to 
methods of numerical modeling of
families of stochastic integrals like (\ref{ito}), (\ref{str}), which do not
use partitioning of the aforementioned
interval of integration $[t, T]$ and converge in the mean-square sense. 
It was suggested in \cite{Mi2} to use converging 
in the mean-square sense trigonometric Fourier 
expansions of the Wiener processes, 
which underlie the iterated stochastic integral.
By this method, 
the mean-square approximations of
the simplest integrals like (\ref{ito}) of multiplicities 1 and 2
($k=2;$ $\psi_1(s),$ $\psi_2(s)\equiv 1;$ $i_1, i_2=0, 1,\ldots,m$)
were obtained in \cite{Mi2}. These approximations were used in \cite{Mi2} 
to construct a numerical method for
the Ito SDE (\ref{1.5.2}), which under certain conditions \cite{Mi2}
has the 
order 1.0 of the mean-square convergence 
and is known as the Milstein method.

A more general method of the mean-square approximation of 
the stochastic integrals like (\ref{str}), which
based on the generalized iterated Fourier series was proposed in 
\cite{3}, \cite{4}. 
It enables one to use the complete orthonormal
systems of Legendre polynomials and trigonometric functions
in the space $L_2([t, T])$. 
In virtue of its characteristics, the method from \cite{Mi2} admits 
the application of only trigonometric basis
functions.

In \cite{KlPl2}, \cite{KPS}, \cite{KPW}, \cite{Zapad-9}
an attempt was made to extend the method 
from \cite{Mi2} to the stochastic integrals
like (\ref{str}) for 
$k=3;$ $\psi_1(s),\ldots,\psi_3(s)\equiv 1;$ $i_1,\ldots,i_3=0, 1,\ldots,m.$

We note that the methods \cite{KlPl2}, \cite{KPS}, \cite{KPW}, 
\cite{Zapad-9}
($k=3$) and \cite{3} ($k\ge 3$)
lead to iterated application of the operation of limit transition.
As a result, these methods allow us to represent the integrals 
(\ref{str})
as iterated series of products of standard Gaussian
random variables (the operation of passing to the limit is 
carried out iteratively). 
This fact is
essential and imposes some constraints related with the 
method of summation of the aforementioned
series \cite{KlPl2}, \cite{KPS}, \cite{3},
\cite{KPW}, \cite{Zapad-9} if we consider the 
stochastic integrals like (\ref{ito}), (\ref{str}) of
multiplicities 3 and higher 
(we mean here at least triple 
integration over the Wiener processes).
Additionally, the aforementioned methods in virtue of their 
features prevent precise calculation
of the mean-square error of approximation with the exception 
of the simplest iterated stochastic
integrals of multip\-li\-ci\-ty 2. This means that at the stage 
of realization of the numerical methods
for Ito SDEs, possibly, one will need to allow for 
the redundant terms of the expansions of
iterated stochastic integrals, which increases 
the computing costs and reduces efficiency of the
numerical methods.

We notice \cite{KlPl2}, \cite{7} that to construct numerical methods for the 
Ito SDE (\ref{1.5.2}) having orders 1.5 and 2.0 of
strong convergence one has to approximate 
(proceeding from the mean-square convergence
criterion) the stochastic integrals not only of multiplicities 
1 and 2, but also 3 and 4 from the
families (\ref{ito}), (\ref{str}). Some 
publications \cite{KlPl2}, \cite{Mi2},
\cite{Mi3} contain the 
aforementioned numerical schemes with
orders 1.5 and 2.0 of strong convergence but without the 
contained in them efficient procedures
of the mean-square approximation of iterated stochastic integrals 
for the case of a multidimensional
Wiener process, which corresponds to 
$i_1,\ldots,i_4=1,\ldots,m$
in (\ref{ito}), (\ref{str}). Part of publications
(see, for example, \cite{KlPl2}, \cite{Mi3})
contain representations of the stochastic integrals of 
multiplicities 3 and 4 like (\ref{ito}), (\ref{str})
only for the simplest case $\psi_1(s),\ldots,\psi_4(s)\equiv 1,$
$i_1=\ldots=i_4$
(representations based 
on the Hermit polynomials). Some publications
\cite{Mi3} use other simplifying assumptions about the Ito SDE (\ref{1.5.2}). 
For example, assumptions
are made about additivity of the stochastic perturbation or its 
smallness, which corresponds,
respectively, to 
$B({\bf x}, t)\equiv C(t)$ or
$B({\bf x}, t)\equiv \varepsilon D({\bf x}, t).$
Here, $\varepsilon>0$ is a fixed small number and
$C: [0, T]\to\mathbb{R}^{n\times m}$, 
$D: \mathbb{R}^n\times[0, T]\to\mathbb{R}^{n\times m}.$

In the case at hand, 
the problem of efficient joint
numerical modeling of the iterated stochastic integrals from the 
families (\ref{ito}), (\ref{str}) becomes 
simpler due to the absence of some terms in the expressions of 
the numerical methods or
the possibility of disregarding some of the aforementioned terms. 
Also, one may encounter approximation 
method \cite{Al} for iterated stochastic integrals of 
mul\-tip\-li\-ci\-ty 3 from the familiy (\ref{ito}) 
for $\psi_1(s),$ $\psi_2(s),$ $\psi_3(s)$ $\equiv 1$ 
$(i_1,i_2,i_3=1,\ldots,m)$ based on partitioning of 
the interval of integration $[t, T]$ of the 
iterated stochastic integrals and using 
multiple integral sums whose disadvantages were mentioned above.

The present paper is devoted to the development of efficient 
procedures 
for joint
numerical modeling of the iterated stochastic integrals from 
the families (\ref{ito}), (\ref{str}) in accordance with the 
mean-square criterion of convergence. At that we do not use 
any essential simplifying assumptions, that is, the
Wiener process involved in the Ito SDE (\ref{1.5.2}) is assumed to be 
the multidimensional one which corresponds to
the condition $i_1,\ldots,i_k=0, 1,\ldots,m$ in (\ref{ito}), (\ref{str}).
In addition, it is assumed that the stochastic perturbation is
nonadditive 
(the simplifying assumptions about the function 
$B: \mathbb{R}^n\times[0, T]\to\mathbb{R}^{n\times m}$ 
involved in (\ref{1.5.2}) are not
introduced).
Additionally, 
the functions $\psi_1(s),\ldots,\psi_k(s)$ in (\ref{ito}), (\ref{str}) are, 
generally speaking,
assumed to be different.
Moreover, the assumption of commutativity \cite{KlPl2}, \cite{KPS} of
the stochastic perturbation is also not introduced.

More precisely, in this paper we consider the method 
of the mean-square approximation of iterated Ito stochastic 
integrals 
from the family (\ref{ito}), which is based
on the generalized multiple
(not iterated) Fourier series converging in the sense of norm
in Hilbert space $L_2([t, T]^k)$ $(k\in\mathbb{N})$ \cite{7} (2006), 
\cite{8}-\cite{31}.
Multiple series (the operation of limit transition 
is implemented only once) are more convenient 
for approximation than the iterated ones
(iterated application of the operation of limit
transition), 
since partial sums of multiple series converge for any possible case of  
convergence to infinity of their upper limits of summation 
(let us denote them as $p_1,\ldots, p_k$). 
For example, when $p_1=\ldots=p_k=p\to\infty$. 
For iterated series, the condition $p_1=\ldots=p_k=p\to\infty$ obviously 
does not guarantee the convergence of this series.
However, 
in \cite{KlPl2}
(Sect.~5.8, pp.~202--204), \cite{KPS} (pp.~82-84),
\cite{KPW} (pp.~438-439),  
\cite{Zapad-9} (pp.~263-264) 
the authors use (without rigorous proof)
the condition $p_1=p_2=p_3=p\to\infty$
within the frames of the approach
based on the Karhunen--Loeve expansion of the Brownian bridge
process \cite{Mi2} together with the Wong--Zakai approximation
\cite{W-Z-1}-\cite{Watanabe}. See discussions 
in \cite{20a} (Sect.~2.18, 6.2), 
\cite{20aa}, \cite{20aaa} (Sect.~2.6.2, 6.2)
for details.

\vspace{5mm}

\section{Numerical Schemes With the Orders 1.0, 1.5, and 2.0 of Strong Convergence}

\vspace{5mm}

Consider the partition $\{\tau_j\}_{j=0}^N$ of the segment $[0, T]$ 
with the partition rank $\Delta_N$ such that

$$
0=\tau_0<\tau_1<\ldots <\tau_N=T.
$$ 

\vspace{3mm}

Denote by 
${\bf y}_{\tau_j}\stackrel{\sf def}{=}
{\bf y}_{j};$\ $j=0, 1,\ldots,N$
the discrete approximation of
the process ${\bf x}_t,$ $t\in[0,T]$ (solution of the Ito SDE (\ref{1.5.2}) 
corresponding to the maximal step of
discretization $\Delta_N$.

\vspace{2mm}

{\bf Definition 1}\ \cite{KlPl2}.\ {\it We will say that
the discrete approximation {\rm (}numerical method{\rm )}
${\bf y}_{j};$\ $j=0, 1,\ldots,N$ corresponding
to the maximal step of discretization $\Delta_N$ converges 
strongly with the order $\gamma>0$
at the time instant $T$ to the process ${\bf x}_t,$ $t\in[0,T]$
if there exist
a constant $C>0$ independent of $\Delta_N$
and a number $\delta>0$ such that

\begin{equation}
\label{xz1}
{\sf M}\left\{\left|{\bf x}_T-{\bf y}_T\right|\right\}\le
C(\Delta_N)^{\gamma}
\end{equation}

\vspace{3mm}
\noindent
for all 
$\Delta_N\in(0, \delta).$}

\vspace{2mm}

We note that the authors of some publications \cite{Mi2}, \cite{Mi3} prefer 
to consider the mean-square convergence
instead of the strong convergence.

\vspace{2mm}

{\bf Definition 2}\ \cite{Mi2}, \cite{Mi3}.\ {\it We will say that
the numerical method 
${\bf y}_{j};$\ $j=0, 1,\ldots,N$ converges in the mean-square
sense with the order $\gamma>0$ to the process ${\bf x}_t,$ $t\in[0,T]$ 
if there exist
a constant $C>0$ independent                         
of $\Delta_N,$ $j$ and a number $\delta>0$ such that

$$
\left({\sf M}
\left\{\left|{\bf x}_j-{\bf y}_j\right|^2\right\}\right)^{1/2}\le
C(\Delta_N)^{\gamma}
$$

\vspace{3mm}
\noindent
for all $\Delta_N\in(0, \delta).$}

\vspace{2mm}

Here, ${\bf x}_{\tau_j}\stackrel{\sf def}{=}
{\bf x}_{j};$\ $j=0, 1,\ldots,N.$ 

We notice that sometimes the condition (\ref{xz1})
in Definition 1 is replaced by
the condition \cite{KlPl2}

$$
{\sf M}\left\{\left|{\bf x}_j-{\bf y}_j\right|\right\}\le
C(\Delta_N)^{\gamma}\ \ \ (j=0, 1,\ldots,N)
$$ 

\vspace{3mm}

At that, the constant $C$ is independent
of $\Delta_N$ and $j.$

Strong convergence follows, obviously, from the mean-square convergence 
in virtue of the Lyapunov
inequality. In what follows, we rely on Definition 1 of strong convergence.

Consider the following 
explicit one-step numerical method

$$
{\bf y}_{p+1}={\bf y}_p+\sum_{i=1}^{m}B_{i}
{\hat I}_{(0)\tau_{p+1},\tau_p}^{(i)}+\Delta{\bf a}
+\sum_{i,j=1}^{m}G_j
B_{i}{\hat I}_{(00)\tau_{p+1},\tau_p}^{(ji)}+
$$

\vspace{1mm}
$$
+
\sum_{i=1}^{m}\Biggl(G_i{\bf a}\left(
\Delta {\hat I}_{(0)\tau_{p+1},\tau_p}^{(i)}+
{\hat I}_{(1)\tau_{p+1},\tau_p}^{(i)}\right)
-LB_{i}{\hat I}_{(1)\tau_{p+1},\tau_p}^{(i)}\Biggr)+
$$

\vspace{1mm}
\begin{equation}
\label{4.18}
+\sum_{i,j,l=1}^{m} G_lG_j
B_{i}{\hat I}_{(000)\tau_{p+1},\tau_p}^{(lji)}
+\frac{\Delta^2}{2}L{\bf a}
\end{equation}

\vspace{4mm}
\noindent
corresponding to the 
constant discretization step
$\Delta=T/N$ 
($\tau_p=p\Delta;$\
$p=0,1,\ldots,N;$ $N>1$),
where
${\hat I}_{(l_1\ldots l_k)s,t}^{(i_1\ldots i_k)}$ 
denotes approximation of the iterated Ito stochastic integral 

\vspace{-1mm}
\begin{equation}
\label{ll1}
I_{(l_1\ldots l_k)s,t}^{(i_1\ldots i_k)}=
 \int\limits^ {s} _ {t} (t-\tau _
{k}) ^ {l_ {k}} 
\ldots \int\limits^ {\tau _ {2}} _ {t} (t-\tau _ {1}) ^ {l_ {1}} d
{\bf f} ^ {(i_ {1})} _ {\tau_ {1}} \ldots 
d {\bf f} _ {\tau_ {k}} ^ {(i_ {k})},
\end{equation}

\vspace{3mm}
\noindent
and
$$
L= {\partial \over \partial t}
+ \sum^ {n} _ {i=1} {\bf a}_i ({\bf x},  t) 
{\partial  \over  \partial  {\bf  x}_i}
+ {1\over 2} \sum^ {m} _ {j=1} \sum^ {n} _ {l,i=1}
B_{lj} ({\bf x}, t) B_{ij} ({\bf x}, t) {\partial
^{2} \over \partial {\bf x}_l \partial {\bf x}_i},
$$

\vspace{3mm}
$$
G_i = \sum^ {n} _ {j=1} B_{ji} ({\bf x}, t)
{\partial  \over \partial {\bf x}_j}\ \ \
(i=1,\ldots,m),
$$

\vspace{5mm}
\noindent
$l_1,\ldots, l_k=0, 1, 2\ldots;$\
$i_1,\ldots, i_k=1,\ldots,m;$\ $k=1, 2,\ldots$;\
$B_i$ and $B_{ij}$ are, respectively, the $i$th column
and $ij$th element of the matrix function $B$; ${\bf a}_i$ and ${\bf x}_i$ 
are, 
respectively, the $i$th components of the
vector function ${\bf a}$ and column ${\bf x}$; the columns 

$$
B_{i},\ \ \ {\bf a},\ \ \ G_jB_{i},\ \ \
G_i{\bf a},\ \ \ LB_{i},\ \ \ G_lG_jB_{i},\ \ \ L{\bf a}
$$ 

\vspace{3mm}
\noindent
are calculated at
the point $({\bf y}_p,p).$

The numerical scheme (\ref{4.18}) can be found, for example, in a somewhat 
different form in \cite{KlPl2}, \cite{Mi2}, \cite{Mi3}.
The difference here lies in that the author of this work 
used in (\ref{4.18})
the relation

\vspace{-1mm}
\begin{equation}
\label{rr}
\Delta {I}_{(0)\tau_{p+1},\tau_p}^{(i)}+
{I}_{(1)\tau_{p+1},\tau_p}^{(i)}=
\int\limits_{\tau_p}^{\tau_{p+1}}
\int\limits_{\tau_p}^{\tau}d{\bf f}_s^{(i)}d\tau
\end{equation}

\vspace{3mm}
\noindent
which follows with probability 1 from the Ito formula and 
enables one to reduce by one the number
of iterated Ito stochastic integrals to be approximated. 
This is due to the fact that the Ito stochastic
integral on the right-hand side of (\ref{rr}) is expressed as a linear
combination of
the Ito stochastic integrals 

$$
{I}_{(0)\tau_{p+1},\tau_p}^{(i)}\ \ \ \hbox{and}\ \ \
{I}_{(1)\tau_{p+1},\tau_p}^{(i)},
$$

\vspace{3mm}
\noindent
whose approximations are already included in the right-hand side of (\ref{4.18}).

It is common knowledge that under certain conditions \cite{KlPl2} the 
discrete approximation (numerical
method) (\ref{4.18}) has the order $1.5$ of strong convergence. Among the 
aforementioned conditions we
note only the condition for approximations of the 
iterated Ito stochastic integrals involved
in (\ref{4.18})

\vspace{-1mm}
\begin{equation}
\label{4.3}
{\sf M}\Biggl\{\Biggl(I_{(l_{1}\ldots l_{k})\tau_{p+1},\tau_p}
^{(i_{1}\ldots i_{k})} 
-\hat I_{(l_{1}\ldots l_{k})\tau_{p+1},\tau_p}^{(i_{1}\ldots i_{k})}
\Biggr)^2\Biggr\}\le C\Delta^{r},
\end{equation}

\vspace{3mm}
\noindent
where $r=4$ and the constant $C$ is independent of $\Delta,$ 
because the present paper deals mostly with
the approximation of the aforementioned stochastic integrals.

Conditions somewhat different from \cite{KlPl2} are given 
in \cite{Mi3}. 
Under them the numerical
method (\ref{4.18}) has the order $1.5$ 
of the mean-square convergence. 

Note that the Milstein method \cite{Mi2} (method with the order 1.0
of strong convergence) corresponds to the first line in (\ref{4.18}).

Consider the explicit one-step
numerical method with the order 2.0 of strong convergence given by

\vspace{1mm}
$$
{\bf y}_{p+1}={\bf y}_p+\sum_{i=1}^{m}B_{i}
{\hat I}_{(0)\tau_{p+1},\tau_p}^{(i)}+\Delta{\bf a}
+\sum_{i,j=1}^{m}G_j
B_{i}{\hat I}_{(00)\tau_{p+1},\tau_p}^{(ji)}+
$$

\vspace{1mm}
$$
+
\sum_{i=1}^{m}\Biggl(G_i{\bf a}\left(
\Delta {\hat I}_{(0)\tau_{p+1},\tau_p}^{(i)}+
{\hat I}_{(1)\tau_{p+1},\tau_p}^{(i)}\right)
-LB_{i}{\hat I}_{(1)\tau_{p+1},\tau_p}^{(i)}\Biggr)+
$$

\vspace{1mm}
$$
+\sum_{i,j,l=1}^{m} G_lG_j
B_{i}{\hat I}_{(000)\tau_{p+1},\tau_p}^{(lji)}
+\frac{\Delta^2}{2}L{\bf a}+
$$

$$
+
\sum_{i,j=1}^{m}
\Biggl(G_0^{(j)}LB_{i}\left(
{\hat I}_{(10)\tau_{p+1},\tau_p}^{(ji)}-
{\hat I}_{(01)\tau_{p+1},\tau_p}^{(ji)}
\right)
-LG_j B_{i}{\hat I}_{(10)\tau_{p+1},\tau_p}^{(ji)}
+\Biggr.
$$

\vspace{1mm}
$$
\Biggl.+G_jG_i{\bf a}\left(
{\hat I}_{(01)\tau_{p+1},\tau_p}^{(ji)}+
\Delta {\hat I}_{(00)\tau_{p+1},\tau_p}^{(ji)}
\right)\Biggr)
+
$$

\vspace{1mm}
\begin{equation}
\label{4.35}
+\sum_{i,j,l,r=1}^{m}G_rG_lG_j B_{i}
{\hat I}_{(0000)\tau_{p+1},\tau_p}^{(r l j i)},
\end{equation}

\vspace{7mm}
\noindent
where notation corresponds to (\ref{4.18}).

The numerical scheme (\ref{4.35}) can be found in another representation 
in \cite{KlPl2}, \cite{Mi3}. In this case
the distinctions are due to the fact that along with (\ref{rr}) the author 
used in (\ref{4.35}) the equalities

\begin{equation}
\label{rr5}
{I}_{(01)\tau_{p+1},\tau_p}^{(ji)}+
\Delta{I}_{(00)\tau_{p+1},\tau_p}^{(ji)}=
\int\limits_{\tau_p}^{\tau_{p+1}}
\int\limits_{\tau_p}^{\theta}
\int\limits_{\tau_p}^{\tau}d{\bf f}_s^{(j)}d{\bf f}_{\tau}^{(i)}d\theta
\end{equation}

\begin{equation}
\label{rr6}
{I}_{(10)\tau_{p+1},\tau_p}^{(ji)}-
{I}_{(01)\tau_{p+1},\tau_p}^{(ji)}=
\int\limits_{\tau_p}^{\tau_{p+1}}
\int\limits_{\tau_p}^{\theta}
\int\limits_{\tau_p}^{\tau}d{\bf f}_s^{(j)}d\tau d{\bf f}_{\theta}^{(i)},
\end{equation}

\vspace{3mm}
\noindent
which follow with probability 1 from the Ito formula and enable one 
to reduce by one more unit
the number of iterated Ito stochastic integrals to be approximated. 
This is due to the fact that
the Ito stochastic integrals on the right-hand sides of (\ref{rr5}) and 
(\ref{rr6}) are 
expressed as linear
combinations of
the Ito stochastic integrals 

\vspace{-1mm}
$$
{I}_{(01)\tau_{p+1},\tau_p}^{(ji)},\ \ \
{I}_{(10)\tau_{p+1},\tau_p}^{(ji)},\ \ \
{I}_{(00)\tau_{p+1},\tau_p}^{(ji)},
$$

\vspace{3mm}
\noindent
whose approximations are already included in the right-hand side of 
(\ref{4.35}).

We notice that under certain conditions \cite{KlPl2} the numerical 
method (\ref{4.35}) has the order $2.0$ of strong
convergence. Among the aforementioned conditions we mark only 
the condition (\ref{4.3}) for $r=5$
intended for approximations of the iterated Ito stochastic integrals 
included in (\ref{4.35}).

Some modifications 
of the numerical methods (\ref{4.18})
and (\ref{4.35}) were constructed
in \cite{KlPl2}, \cite{Mi3}. 
Among which there are finite-difference methods
of the Runge--Kutta type as well
as the implicit and two-step methods (also see \cite{7}, \cite{8}-\cite{11},
\cite{19}-\cite{20aaa}). In all aforementioned methods, 
however, a need arises for
efficient joint mean-square approximation of the iterated 
Ito stochastic integrals of multiplicities 1 to 4.
The collection of these integrals is the same as in 
the numerical methods (\ref{4.18}) and (\ref{4.35}).

\vspace{5mm}

\section{Expansion of Iterated Ito Stochastic Integrals
of Multiplicity $k$ $(k\in \mathbb{N})$
Based on Generalized Multiple Fourier Series}

\vspace{5mm}

An efficient mean-square approximation method for the 
iterated Ito stochastic 
integrals like (\ref{ito}) was
proposed and developed by the author of this article in \cite{7}, 
\cite{8}-\cite{31} (see Theorems 1, 2 below). This method based
on the generalized multiple Fourier series
converging in the mean-square sense in the space
$L_2([t, T]^k),$ $k\in\mathbb{N}.$ 
At that the method \cite{7}, 
\cite{8}-\cite{31} allows to use different
complete orthonormal systems of functions
in the space $L_2([t, T]^k),$ $k\in\mathbb{N}.$
In this article, we use the system of Legendre polynomials,
which has a series of advantages over the system of 
trigonometric functions in the framework of the
considered problem \cite{29}, \cite{30}.
Moreover, in this method the passage to the 
limit is carried out only once, which leads
to a correct choice of the lengths of sequences of the standard 
Gaussian random variables required
to approximate the iterated Ito stochastic integrals.

Suppose that every $\psi_l(\tau)$ $(l=1,\ldots,k)$ is a 
nonrandom function from the space $L_2([t, T])$.
Define the following function on the hypercube $[t, T]^k$

\vspace{-1mm}
\begin{equation}
\label{ppp}
K(t_1,\ldots,t_k)=
\begin{cases}
\psi_1(t_1)\ldots \psi_k(t_k)\ &\hbox{for}\ \ t_1<\ldots<t_k\\
~\\
~\\
0\ &\hbox{otherwise}
\end{cases},\ \ \ \ t_1,\ldots,t_k\in[t, T],\ \ \ \ k\ge 2,
\end{equation}

\vspace{5mm}
\noindent
and 
$K(t_1)\equiv\psi_1(t_1)$ for $t_1\in[t, T].$

Suppose that $\{\phi_j(x)\}_{j=0}^{\infty}$
is a complete orthonormal system of functions in the space
$L_2([t, T])$. 
The function $K(t_1,\ldots,t_k)$ belongs to the space
$L_2([t, T]^k).$
At this situation it is well known that the generalized 
multiple Fourier series 
of $K(t_1,\ldots,t_k)\in L_2([t, T]^k)$ is converging 
to $K(t_1,\ldots,t_k)$ in the hypercube $[t, T]^k$ in 
the mean-square sense, i.e.

$$
\hbox{\vtop{\offinterlineskip\halign{
\hfil#\hfil\cr
{\rm lim}\cr
$\stackrel{}{{}_{p_1,\ldots,p_k\to \infty}}$\cr
}} }\Biggl\Vert
K(t_1,\ldots,t_k)-
\sum_{j_1=0}^{p_1}\ldots \sum_{j_k=0}^{p_k}
C_{j_k\ldots j_1}\prod_{l=1}^{k} \phi_{j_l}(t_l)
\Biggr\Vert_{L_2([t,T]^k)}=0,
$$

\vspace{4mm}
\noindent
where
\begin{equation}
\label{ppppa}
C_{j_k\ldots j_1}=\int\limits_{[t,T]^k}
K(t_1,\ldots,t_k)\prod_{l=1}^{k}\phi_{j_l}(t_l)dt_1\ldots dt_k
\end{equation}

\vspace{3mm}
\noindent
is the Fourier coefficient,

$$
\left\Vert f\right\Vert_{L_2([t,T]^k)}=\left(\int\limits_{[t,T]^k}
f^2(t_1,\ldots,t_k)dt_1\ldots dt_k\right)^{1/2},
$$

\vspace{4mm}
\noindent
and the Parceval equality

\begin{equation}
\label{sos2z}
\int\limits_{[t,T]^k}
K^2(t_1,\ldots,t_k)dt_1\ldots dt_k
=\hbox{\vtop{\offinterlineskip\halign{
\hfil#\hfil\cr
{\rm lim}\cr
$\stackrel{}{{}_{p_1,\ldots,p_k\to \infty}}$\cr
}} }\sum_{j_1=0}^{p_1}\ldots\sum_{j_k=0}^{p_k}C_{j_k\ldots j_1}^2
\end{equation}

\vspace{4mm}
\noindent
takes place.

\vspace{4mm}

Consider the partition $\{\tau_j\}_{j=0}^N$ of $[t,T]$ such that

\vspace{1mm}

\begin{equation}
\label{1111}
t=\tau_0<\ldots <\tau_N=T,\ \ \
\Delta_N=
\hbox{\vtop{\offinterlineskip\halign{
\hfil#\hfil\cr
{\rm max}\cr
$\stackrel{}{{}_{0\le j\le N-1}}$\cr
}} }\Delta\tau_j\to 0\ \ \hbox{if}\ \ N\to \infty,\ \ \
\Delta\tau_j=\tau_{j+1}-\tau_j.
\end{equation}

\vspace{4mm}

{\bf Theorem 1}\ \cite{7} (2006), \cite{8}-\cite{31}.
{\it Suppose that
every $\psi_l(\tau)$ $(l=1,\ldots, k)$ is a continuous nonrandom function on 
$[t, T]$ and
$\{\phi_j(x)\}_{j=0}^{\infty}$ is a complete orthonormal system  
of continuous functions in the space $L_2([t,T]).$ Then

\vspace{1mm}
$$
J[\psi^{(k)}]_{T,t}\  =\ 
\hbox{\vtop{\offinterlineskip\halign{
\hfil#\hfil\cr
{\rm l.i.m.}\cr
$\stackrel{}{{}_{p_1,\ldots,p_k\to \infty}}$\cr
}} }\sum_{j_1=0}^{p_1}\ldots\sum_{j_k=0}^{p_k}
C_{j_k\ldots j_1}\Biggl(
\prod_{l=1}^k\zeta_{j_l}^{(i_l)}\ -
\Biggr.
$$

\vspace{3mm}
\begin{equation}
\label{tyyy}
-\ \Biggl.
\hbox{\vtop{\offinterlineskip\halign{
\hfil#\hfil\cr
{\rm l.i.m.}\cr
$\stackrel{}{{}_{N\to \infty}}$\cr
}} }\sum_{(l_1,\ldots,l_k)\in {\rm G}_k}
\phi_{j_{1}}(\tau_{l_1})
\Delta{\bf w}_{\tau_{l_1}}^{(i_1)}\ldots
\phi_{j_{k}}(\tau_{l_k})
\Delta{\bf w}_{\tau_{l_k}}^{(i_k)}\Biggr),
\end{equation}

\vspace{6mm}
\noindent
where $J[\psi^{(k)}]_{T,t}$ is defined by {\rm (\ref{ito}),}

\vspace{2mm}
$$
{\rm G}_k={\rm H}_k\backslash{\rm L}_k,\ \ \
{\rm H}_k=\{(l_1,\ldots,l_k):\ l_1,\ldots,l_k=0,\ 1,\ldots,N-1\},
$$

\vspace{1mm}
$$
{\rm L}_k=\{(l_1,\ldots,l_k):\ l_1,\ldots,l_k=0,\ 1,\ldots,N-1;\
l_g\ne l_r\ (g\ne r);\ g, r=1,\ldots,k\},
$$

\vspace{6mm}
\noindent
${\rm l.i.m.}$ is a limit in the mean-square sense$,$
$i_1,\ldots,i_k=0,1,\ldots,m,$

\begin{equation}
\label{rr23}
\zeta_{j}^{(i)}=
\int\limits_t^T \phi_{j}(s) d{\bf w}_s^{(i)}
\end{equation} 

\vspace{3mm}
\noindent
are independent standard Gaussian random variables
for various
$i$ or $j$ {\rm(}if $i\ne 0${\rm),}
$C_{j_k\ldots j_1}$ is the Fourier coefficient {\rm(\ref{ppppa}),}
$\Delta{\bf w}_{\tau_{j}}^{(i)}=
{\bf w}_{\tau_{j+1}}^{(i)}-{\bf w}_{\tau_{j}}^{(i)}$
$(i=0, 1,\ldots,m),$
$\left\{\tau_{j}\right\}_{j=0}^{N}$ is a partition of
the interval $[t, T],$ which satisfies the condition {\rm (\ref{1111})}.
}

\vspace{4mm}

It was shown in \cite{9}-\cite{16}, \cite{19}-\cite{20aaa}
that 
Theorem 1 is valid for convergence 
in the mean of degree $2n$ ($n\in \mathbb{N}$).
Moreover, the convergence with probability 1
in Theorem 1 is proved in \cite{20a}-\cite{20aaa}, \cite{OK1000}.
In addition, the complete orthonormal systems of Haar and 
Rademacher--Walsh functions in $L_2([t,T])$ also
can be applied in Theorem 1
\cite{7}, \cite{8}-\cite{16}, \cite{19}-\cite{20aaa}.
The modification of Theorem 1 for 
complete orthonormal with weigth $r(x)\ge 0$ systems
of functions in the space $L_2([t,T])$ can be found in 
\cite{20}, \cite{20a}-\cite{20aaa}, \cite{26b}. Application of Theorem 1 
and Theorem 2 (see below)
to the approximation of iterated stochastic integrals
with respect to the infinite-dimensional
$Q$-Wiener process contains in \cite{20a}-\cite{20aaa} (Chapter 7),
\cite{31a}, \cite{31}, \cite{OK}, \cite{Kuzh-1}.

In order to evaluate the significance of Theorem 1 for practice we will
demonstrate its transformed particular cases for 
$k=1,\ldots,5$ \cite{7} (2006), 
\cite{8}-\cite{31} (the cases $k=6, 7$ and $k>7$ $(k\in \mathbb{N})$ 
can also be found in these papers)

\begin{equation}
\label{a1}
J[\psi^{(1)}]_{T,t}
=\hbox{\vtop{\offinterlineskip\halign{
\hfil#\hfil\cr
{\rm l.i.m.}\cr
$\stackrel{}{{}_{p_1\to \infty}}$\cr
}} }\sum_{j_1=0}^{p_1}
C_{j_1}\zeta_{j_1}^{(i_1)},
\end{equation}

\vspace{1mm}
\begin{equation}
\label{a2}
J[\psi^{(2)}]_{T,t}
=\hbox{\vtop{\offinterlineskip\halign{
\hfil#\hfil\cr
{\rm l.i.m.}\cr
$\stackrel{}{{}_{p_1,p_2\to \infty}}$\cr
}} }\sum_{j_1=0}^{p_1}\sum_{j_2=0}^{p_2}
C_{j_2j_1}\Biggl(\zeta_{j_1}^{(i_1)}\zeta_{j_2}^{(i_2)}
-{\bf 1}_{\{i_1=i_2\ne 0\}}
{\bf 1}_{\{j_1=j_2\}}\Biggr),
\end{equation}

\vspace{5mm}

$$
J[\psi^{(3)}]_{T,t}=
\hbox{\vtop{\offinterlineskip\halign{
\hfil#\hfil\cr
{\rm l.i.m.}\cr
$\stackrel{}{{}_{p_1,p_2,p_3\to \infty}}$\cr
}} }\sum_{j_1=0}^{p_1}\sum_{j_2=0}^{p_2}\sum_{j_3=0}^{p_3}
C_{j_3j_2j_1}\Biggl(
\zeta_{j_1}^{(i_1)}\zeta_{j_2}^{(i_2)}\zeta_{j_3}^{(i_3)}
-\Biggr.
$$

\begin{equation}
\label{a3}
-\Biggl.
{\bf 1}_{\{i_1=i_2\ne 0\}}
{\bf 1}_{\{j_1=j_2\}}
\zeta_{j_3}^{(i_3)}
-{\bf 1}_{\{i_2=i_3\ne 0\}}
{\bf 1}_{\{j_2=j_3\}}
\zeta_{j_1}^{(i_1)}-
{\bf 1}_{\{i_1=i_3\ne 0\}}
{\bf 1}_{\{j_1=j_3\}}
\zeta_{j_2}^{(i_2)}\Biggr),
\end{equation}

\vspace{8mm}

$$
J[\psi^{(4)}]_{T,t}
=
\hbox{\vtop{\offinterlineskip\halign{
\hfil#\hfil\cr
{\rm l.i.m.}\cr
$\stackrel{}{{}_{p_1,\ldots,p_4\to \infty}}$\cr
}} }\sum_{j_1=0}^{p_1}\ldots \sum_{j_4=0}^{p_4}
C_{j_4\ldots j_1}\Biggl(
\prod_{l=1}^4\zeta_{j_l}^{(i_l)}
\Biggr.
-
$$
$$
-
{\bf 1}_{\{i_1=i_2\ne 0\}}
{\bf 1}_{\{j_1=j_2\}}
\zeta_{j_3}^{(i_3)}
\zeta_{j_4}^{(i_4)}
-
{\bf 1}_{\{i_1=i_3\ne 0\}}
{\bf 1}_{\{j_1=j_3\}}
\zeta_{j_2}^{(i_2)}
\zeta_{j_4}^{(i_4)}-
$$
$$
-
{\bf 1}_{\{i_1=i_4\ne 0\}}
{\bf 1}_{\{j_1=j_4\}}
\zeta_{j_2}^{(i_2)}
\zeta_{j_3}^{(i_3)}
-
{\bf 1}_{\{i_2=i_3\ne 0\}}
{\bf 1}_{\{j_2=j_3\}}
\zeta_{j_1}^{(i_1)}
\zeta_{j_4}^{(i_4)}-
$$
$$
-
{\bf 1}_{\{i_2=i_4\ne 0\}}
{\bf 1}_{\{j_2=j_4\}}
\zeta_{j_1}^{(i_1)}
\zeta_{j_3}^{(i_3)}
-
{\bf 1}_{\{i_3=i_4\ne 0\}}
{\bf 1}_{\{j_3=j_4\}}
\zeta_{j_1}^{(i_1)}
\zeta_{j_2}^{(i_2)}+
$$
$$
+
{\bf 1}_{\{i_1=i_2\ne 0\}}
{\bf 1}_{\{j_1=j_2\}}
{\bf 1}_{\{i_3=i_4\ne 0\}}
{\bf 1}_{\{j_3=j_4\}}
+
{\bf 1}_{\{i_1=i_3\ne 0\}}
{\bf 1}_{\{j_1=j_3\}}
{\bf 1}_{\{i_2=i_4\ne 0\}}
{\bf 1}_{\{j_2=j_4\}}+
$$
\begin{equation}
\label{a4}
+\Biggl.
{\bf 1}_{\{i_1=i_4\ne 0\}}
{\bf 1}_{\{j_1=j_4\}}
{\bf 1}_{\{i_2=i_3\ne 0\}}
{\bf 1}_{\{j_2=j_3\}}\Biggr),
\end{equation}

\vspace{8mm}

$$
J[\psi^{(5)}]_{T,t}
=\hbox{\vtop{\offinterlineskip\halign{
\hfil#\hfil\cr
{\rm l.i.m.}\cr
$\stackrel{}{{}_{p_1,\ldots,p_5\to \infty}}$\cr
}} }\sum_{j_1=0}^{p_1}\ldots\sum_{j_5=0}^{p_5}
C_{j_5\ldots j_1}\Biggl(
\prod_{l=1}^5\zeta_{j_l}^{(i_l)}
-\Biggr.
$$
$$
-
{\bf 1}_{\{i_1=i_2\ne 0\}}
{\bf 1}_{\{j_1=j_2\}}
\zeta_{j_3}^{(i_3)}
\zeta_{j_4}^{(i_4)}
\zeta_{j_5}^{(i_5)}-
{\bf 1}_{\{i_1=i_3\ne 0\}}
{\bf 1}_{\{j_1=j_3\}}
\zeta_{j_2}^{(i_2)}
\zeta_{j_4}^{(i_4)}
\zeta_{j_5}^{(i_5)}-
$$
$$
-
{\bf 1}_{\{i_1=i_4\ne 0\}}
{\bf 1}_{\{j_1=j_4\}}
\zeta_{j_2}^{(i_2)}
\zeta_{j_3}^{(i_3)}
\zeta_{j_5}^{(i_5)}-
{\bf 1}_{\{i_1=i_5\ne 0\}}
{\bf 1}_{\{j_1=j_5\}}
\zeta_{j_2}^{(i_2)}
\zeta_{j_3}^{(i_3)}
\zeta_{j_4}^{(i_4)}-
$$
$$
-
{\bf 1}_{\{i_2=i_3\ne 0\}}
{\bf 1}_{\{j_2=j_3\}}
\zeta_{j_1}^{(i_1)}
\zeta_{j_4}^{(i_4)}
\zeta_{j_5}^{(i_5)}-
{\bf 1}_{\{i_2=i_4\ne 0\}}
{\bf 1}_{\{j_2=j_4\}}
\zeta_{j_1}^{(i_1)}
\zeta_{j_3}^{(i_3)}
\zeta_{j_5}^{(i_5)}-
$$
$$
-
{\bf 1}_{\{i_2=i_5\ne 0\}}
{\bf 1}_{\{j_2=j_5\}}
\zeta_{j_1}^{(i_1)}
\zeta_{j_3}^{(i_3)}
\zeta_{j_4}^{(i_4)}
-{\bf 1}_{\{i_3=i_4\ne 0\}}
{\bf 1}_{\{j_3=j_4\}}
\zeta_{j_1}^{(i_1)}
\zeta_{j_2}^{(i_2)}
\zeta_{j_5}^{(i_5)}-
$$
$$
-
{\bf 1}_{\{i_3=i_5\ne 0\}}
{\bf 1}_{\{j_3=j_5\}}
\zeta_{j_1}^{(i_1)}
\zeta_{j_2}^{(i_2)}
\zeta_{j_4}^{(i_4)}
-{\bf 1}_{\{i_4=i_5\ne 0\}}
{\bf 1}_{\{j_4=j_5\}}
\zeta_{j_1}^{(i_1)}
\zeta_{j_2}^{(i_2)}
\zeta_{j_3}^{(i_3)}+
$$
$$
+
{\bf 1}_{\{i_1=i_2\ne 0\}}
{\bf 1}_{\{j_1=j_2\}}
{\bf 1}_{\{i_3=i_4\ne 0\}}
{\bf 1}_{\{j_3=j_4\}}\zeta_{j_5}^{(i_5)}+
{\bf 1}_{\{i_1=i_2\ne 0\}}
{\bf 1}_{\{j_1=j_2\}}
{\bf 1}_{\{i_3=i_5\ne 0\}}
{\bf 1}_{\{j_3=j_5\}}\zeta_{j_4}^{(i_4)}+
$$
$$
+
{\bf 1}_{\{i_1=i_2\ne 0\}}
{\bf 1}_{\{j_1=j_2\}}
{\bf 1}_{\{i_4=i_5\ne 0\}}
{\bf 1}_{\{j_4=j_5\}}\zeta_{j_3}^{(i_3)}+
{\bf 1}_{\{i_1=i_3\ne 0\}}
{\bf 1}_{\{j_1=j_3\}}
{\bf 1}_{\{i_2=i_4\ne 0\}}
{\bf 1}_{\{j_2=j_4\}}\zeta_{j_5}^{(i_5)}+
$$
$$
+
{\bf 1}_{\{i_1=i_3\ne 0\}}
{\bf 1}_{\{j_1=j_3\}}
{\bf 1}_{\{i_2=i_5\ne 0\}}
{\bf 1}_{\{j_2=j_5\}}\zeta_{j_4}^{(i_4)}+
{\bf 1}_{\{i_1=i_3\ne 0\}}
{\bf 1}_{\{j_1=j_3\}}
{\bf 1}_{\{i_4=i_5\ne 0\}}
{\bf 1}_{\{j_4=j_5\}}\zeta_{j_2}^{(i_2)}+
$$
$$
+
{\bf 1}_{\{i_1=i_4\ne 0\}}
{\bf 1}_{\{j_1=j_4\}}
{\bf 1}_{\{i_2=i_3\ne 0\}}
{\bf 1}_{\{j_2=j_3\}}\zeta_{j_5}^{(i_5)}+
{\bf 1}_{\{i_1=i_4\ne 0\}}
{\bf 1}_{\{j_1=j_4\}}
{\bf 1}_{\{i_2=i_5\ne 0\}}
{\bf 1}_{\{j_2=j_5\}}\zeta_{j_3}^{(i_3)}+
$$
$$
+
{\bf 1}_{\{i_1=i_4\ne 0\}}
{\bf 1}_{\{j_1=j_4\}}
{\bf 1}_{\{i_3=i_5\ne 0\}}
{\bf 1}_{\{j_3=j_5\}}\zeta_{j_2}^{(i_2)}+
{\bf 1}_{\{i_1=i_5\ne 0\}}
{\bf 1}_{\{j_1=j_5\}}
{\bf 1}_{\{i_2=i_3\ne 0\}}
{\bf 1}_{\{j_2=j_3\}}\zeta_{j_4}^{(i_4)}+
$$
$$
+
{\bf 1}_{\{i_1=i_5\ne 0\}}
{\bf 1}_{\{j_1=j_5\}}
{\bf 1}_{\{i_2=i_4\ne 0\}}
{\bf 1}_{\{j_2=j_4\}}\zeta_{j_3}^{(i_3)}+
{\bf 1}_{\{i_1=i_5\ne 0\}}
{\bf 1}_{\{j_1=j_5\}}
{\bf 1}_{\{i_3=i_4\ne 0\}}
{\bf 1}_{\{j_3=j_4\}}\zeta_{j_2}^{(i_2)}+
$$
$$
+
{\bf 1}_{\{i_2=i_3\ne 0\}}
{\bf 1}_{\{j_2=j_3\}}
{\bf 1}_{\{i_4=i_5\ne 0\}}
{\bf 1}_{\{j_4=j_5\}}\zeta_{j_1}^{(i_1)}+
{\bf 1}_{\{i_2=i_4\ne 0\}}
{\bf 1}_{\{j_2=j_4\}}
{\bf 1}_{\{i_3=i_5\ne 0\}}
{\bf 1}_{\{j_3=j_5\}}\zeta_{j_1}^{(i_1)}+
$$
\begin{equation}
\label{a5}
+\Biggl.
{\bf 1}_{\{i_2=i_5\ne 0\}}
{\bf 1}_{\{j_2=j_5\}}
{\bf 1}_{\{i_3=i_4\ne 0\}}
{\bf 1}_{\{j_3=j_4\}}\zeta_{j_1}^{(i_1)}\Biggr),
\end{equation}

\vspace{8mm}
\noindent
where ${\bf 1}_A$ is the indicator of the set $A$.

For further consideration, let us 
consider the generalization of formulas (\ref{a1})--(\ref{a5})                 
for the case of an arbitrary multiplicity $k$ $(k\in\mathbb{N})$ of 
the iterated Ito stochastic integral $J[\psi^{(k)}]_{T,t}$ defined by (\ref{ito}).
In order to do this, let us
introduce some notations. 
Consider the unordered
set $\{1, 2, \ldots, k\}$ 
and separate it into two parts:
the first part consists of $r$ unordered 
pairs (sequence order of these pairs is also unimportant) and the 
second one consists of the 
remaining $k-2r$ numbers.
So, we have

\begin{equation}
\label{leto5007}
(\{
\underbrace{\{g_1, g_2\}, \ldots, 
\{g_{2r-1}, g_{2r}\}}_{\small{\hbox{part 1}}}
\},
\{\underbrace{q_1, \ldots, q_{k-2r}}_{\small{\hbox{part 2}}}
\}),
\end{equation}

\vspace{4mm}
\noindent
where 

\vspace{-2mm}
$$
\{g_1, g_2, \ldots, 
g_{2r-1}, g_{2r}, q_1, \ldots, q_{k-2r}\}=\{1, 2, \ldots, k\},
$$

\vspace{4mm}
\noindent
braces   
mean an unordered 
set, and pa\-ren\-the\-ses mean an ordered set.

We will say that (\ref{leto5007}) is a partition 
and consider the sum with respect to all possible
partitions

\begin{equation}
\label{leto5008}
\sum_{\stackrel{(\{\{g_1, g_2\}, \ldots, 
\{g_{2r-1}, g_{2r}\}\}, \{q_1, \ldots, q_{k-2r}\})}
{{}_{\{g_1, g_2, \ldots, 
g_{2r-1}, g_{2r}, q_1, \ldots, q_{k-2r}\}=\{1, 2, \ldots, k\}}}}
a_{g_1 g_2, \ldots, 
g_{2r-1} g_{2r}, q_1 \ldots q_{k-2r}}.
\end{equation}

\vspace{4mm}

Below there are several examples of sums in the form (\ref{leto5008})

\vspace{2mm}
$$
\sum_{\stackrel{(\{g_1, g_2\})}{{}_{\{g_1, g_2\}=\{1, 2\}}}}
a_{g_1 g_2}=a_{12},
$$

\vspace{3mm}
$$
\sum_{\stackrel{(\{\{g_1, g_2\}, \{g_3, g_4\}\})}
{{}_{\{g_1, g_2, g_3, g_4\}=\{1, 2, 3, 4\}}}}
a_{g_1 g_2 g_3 g_4}=a_{1234} + a_{1324} + a_{2314},
$$

\vspace{3mm}
$$
\sum_{\stackrel{(\{g_1, g_2\}, \{q_1, q_{2}\})}
{{}_{\{g_1, g_2, q_1, q_{2}\}=\{1, 2, 3, 4\}}}}
a_{g_1 g_2, q_1 q_{2}}=
$$

$$
=a_{12,34}+a_{13,24}+a_{14,23}
+a_{23,14}+a_{24,13}+a_{34,12},
$$

\vspace{3mm}
$$
\sum_{\stackrel{(\{g_1, g_2\}, \{q_1, q_{2}, q_3\})}
{{}_{\{g_1, g_2, q_1, q_{2}, q_3\}=\{1, 2, 3, 4, 5\}}}}
a_{g_1 g_2, q_1 q_{2}q_3}
=
$$

$$
=a_{12,345}+a_{13,245}+a_{14,235}
+a_{15,234}+a_{23,145}+a_{24,135}+
$$
$$
+a_{25,134}+a_{34,125}+a_{35,124}+a_{45,123},
$$

\vspace{4mm}
$$
\sum_{\stackrel{(\{\{g_1, g_2\}, \{g_3, g_{4}\}\}, \{q_1\})}
{{}_{\{g_1, g_2, g_3, g_{4}, q_1\}=\{1, 2, 3, 4, 5\}}}}
a_{g_1 g_2, g_3 g_{4},q_1}
=
$$

$$
=
a_{12,34,5}+a_{13,24,5}+a_{14,23,5}+
a_{12,35,4}+a_{13,25,4}+a_{15,23,4}+
$$
$$
+a_{12,54,3}+a_{15,24,3}+a_{14,25,3}+a_{15,34,2}+a_{13,54,2}+a_{14,53,2}+
$$
$$
+
a_{52,34,1}+a_{53,24,1}+a_{54,23,1}.
$$

\vspace{5mm}

Now we can write (\ref{tyyy}) as

\vspace{1mm}

$$
J[\psi^{(k)}]_{T,t}=
\hbox{\vtop{\offinterlineskip\halign{
\hfil#\hfil\cr
{\rm l.i.m.}\cr
$\stackrel{}{{}_{p_1,\ldots,p_k\to \infty}}$\cr
}} }
\sum\limits_{j_1=0}^{p_1}\ldots
\sum\limits_{j_k=0}^{p_k}
C_{j_k\ldots j_1}\Biggl(
\prod_{l=1}^k\zeta_{j_l}^{(i_l)}+\sum\limits_{r=1}^{[k/2]}
(-1)^r \times
\Biggr.
$$

\vspace{3mm}
\begin{equation}
\label{leto6000hh}
\times
\sum_{\stackrel{(\{\{g_1, g_2\}, \ldots, 
\{g_{2r-1}, g_{2r}\}\}, \{q_1, \ldots, q_{k-2r}\})}
{{}_{\{g_1, g_2, \ldots, 
g_{2r-1}, g_{2r}, q_1, \ldots, q_{k-2r}\}=\{1, 2, \ldots, k\}}}}
\prod\limits_{s=1}^r
{\bf 1}_{\{i_{g_{{}_{2s-1}}}=~i_{g_{{}_{2s}}}\ne 0\}}
\Biggl.{\bf 1}_{\{j_{g_{{}_{2s-1}}}=~j_{g_{{}_{2s}}}\}}
\prod_{l=1}^{k-2r}\zeta_{j_{q_l}}^{(i_{q_l})}\Biggr),
\end{equation}

\vspace{5mm}
\noindent
where $[x]$ is an integer part of a real number $x;$
another notations are the same as in Theorem {\bf 1}.

\vspace{2mm}

In particular, from (\ref{leto6000hh}) for $k=5$ we obtain

\vspace{3mm}

$$
J[\psi^{(5)}]_{T,t}=
\hbox{\vtop{\offinterlineskip\halign{
\hfil#\hfil\cr
{\rm l.i.m.}\cr
$\stackrel{}{{}_{p_1,\ldots,p_5\to \infty}}$\cr
}} }\sum_{j_1=0}^{p_1}\ldots\sum_{j_5=0}^{p_5}
C_{j_5\ldots j_1}\Biggl(
\prod_{l=1}^5\zeta_{j_l}^{(i_l)}-\Biggr.
$$

\vspace{2mm}
$$
-
\sum\limits_{\stackrel{(\{g_1, g_2\}, \{q_1, q_{2}, q_3\})}
{{}_{\{g_1, g_2, q_{1}, q_{2}, q_3\}=\{1, 2, 3, 4, 5\}}}}
{\bf 1}_{\{i_{g_{{}_{1}}}=~i_{g_{{}_{2}}}\ne 0\}}
{\bf 1}_{\{j_{g_{{}_{1}}}=~j_{g_{{}_{2}}}\}}
\prod_{l=1}^{3}\zeta_{j_{q_l}}^{(i_{q_l})}+
$$

\vspace{2mm}
$$
+
\sum_{\stackrel{(\{\{g_1, g_2\}, 
\{g_{3}, g_{4}\}\}, \{q_1\})}
{{}_{\{g_1, g_2, g_{3}, g_{4}, q_1\}=\{1, 2, 3, 4, 5\}}}}
{\bf 1}_{\{i_{g_{{}_{1}}}=~i_{g_{{}_{2}}}\ne 0\}}
{\bf 1}_{\{j_{g_{{}_{1}}}=~j_{g_{{}_{2}}}\}}
\Biggl.{\bf 1}_{\{i_{g_{{}_{3}}}=~i_{g_{{}_{4}}}\ne 0\}}
{\bf 1}_{\{j_{g_{{}_{3}}}=~j_{g_{{}_{4}}}\}}
\zeta_{j_{q_1}}^{(i_{q_1})}\Biggr).
$$

\vspace{7mm}
\noindent
The last equality obviously agrees with
(\ref{a5}).

Let us consider the generalization of Theorem 1 for the case
of an arbitrary complete orthonormal systems  
of functions in the space $L_2([t,T])$ 
and $\psi_1(\tau),\ldots,\psi_k(\tau)\in L_2([t, T]).$

\vspace{2mm}

{\bf Theorem~2}\ \cite{20a} (Sect.~1.11), \cite{26a} (Sect.~15).
{\it Suppose that
$\psi_1(\tau),\ldots,\psi_k(\tau)\in L_2([t, T])$ and
$\{\phi_j(x)\}_{j=0}^{\infty}$ is an arbitrary complete orthonormal system  
of functions in the space $L_2([t,T]).$
Then the following expansion

\vspace{1mm}
$$
J[\psi^{(k)}]_{T,t}=
\hbox{\vtop{\offinterlineskip\halign{
\hfil#\hfil\cr
{\rm l.i.m.}\cr
$\stackrel{}{{}_{p_1,\ldots,p_k\to \infty}}$\cr
}} }
\sum\limits_{j_1=0}^{p_1}\ldots
\sum\limits_{j_k=0}^{p_k}
C_{j_k\ldots j_1}\Biggl(
\prod_{l=1}^k\zeta_{j_l}^{(i_l)}+\sum\limits_{r=1}^{[k/2]}
(-1)^r \times
\Biggr.
$$

\vspace{2mm}
\begin{equation}
\label{leto6000}
\times
\sum_{\stackrel{(\{\{g_1, g_2\}, \ldots, 
\{g_{2r-1}, g_{2r}\}\}, \{q_1, \ldots, q_{k-2r}\})}
{{}_{\{g_1, g_2, \ldots, 
g_{2r-1}, g_{2r}, q_1, \ldots, q_{k-2r}\}=\{1, 2, \ldots, k\}}}}
\prod\limits_{s=1}^r
{\bf 1}_{\{i_{g_{{}_{2s-1}}}=~i_{g_{{}_{2s}}}\ne 0\}}
\Biggl.{\bf 1}_{\{j_{g_{{}_{2s-1}}}=~j_{g_{{}_{2s}}}\}}
\prod_{l=1}^{k-2r}\zeta_{j_{q_l}}^{(i_{q_l})}\Biggr)
\end{equation}

\vspace{6mm}
\noindent
con\-verg\-ing in the mean-square sense is valid,
where $[x]$ is an integer part of a real number $x;$
another notations are the same as in Theorem~{\rm 1}.}

\vspace{2mm}

It should be noted that an analogue of Theorem 2 was considered 
in \cite{Rybakov1000}. 
Note that we use another notations 
\cite{20a} (Sect.~1.11), \cite{26a} (Sect.~15)
in comparison with \cite{Rybakov1000}.
Moreover, the proof of an analogue of Theorem 2
from \cite{Rybakov1000} is somewhat different from the proof given in 
\cite{20a} (Sect.~1.11), \cite{26a} (Sect.~15).

\vspace{5mm}

\section{Expansion of Iterated Stratonovich Stochastic
Integrals of Multiplicities 1 to 6}

\vspace{5mm}

As it turned out \cite{12}-\cite{16},
\cite{19}-\cite{20aaa}, \cite{arxiv-5}, \cite{new-art-1-xxy}, \cite{new-art-1xxys}
Theorems 1, 2 can be adapted for the 
iterated Stratonovich stochastic integrals 
(\ref{str}). At that
the expansions of the integrals (\ref{str}) 
turn out to be much simpler than the expansions of 
the iterated Ito stochastic 
integrals (\ref{ito}). 
Let us first present some old results as the following theorem.

\vspace{2mm}

{\bf Theorem 3}\ \cite{12}-\cite{16},
\cite{19}-\cite{20aaa}, \cite{arxiv-5}. 
{\it Assume that the following conditions are fulfilled{\rm :}

{\rm 1}. $\{\phi_j(x)\}_{j=0}^{\infty}$ is the complete orthonormal
system of Legendre polynomials or 
trigonometric functions in the space $L_2([t, T]).$

{\rm 2}. The function $\psi_2(\tau)$ 
is continuously differentiable at the interval $[t, T]$, and 
the functions $\psi_1(\tau),$ $\psi_3(\tau)$ are twice continuously 
differentiable at the interval $[t, T]$ 
{\rm (}in
{\rm (\ref{111})} and {\rm (\ref{feto19000a}))}.

Then, the iterated Stratonovich stochastic
integrals {\rm (\ref{str})} of multiplicities {\rm 2}--{\rm 4} 
are expanded
into the mean-square converging multiple series

\vspace{-1mm}
\begin{equation}
\label{111}
J^{*}[\psi^{(2)}]_{T,t}=
\hbox{\vtop{\offinterlineskip\halign{
\hfil#\hfil\cr
{\rm l.i.m.}\cr
$\stackrel{}{{}_{q_1,q_2\to \infty}}$\cr
}} }\sum_{j_1=0}^{q_1}\sum_{j_2=0}^{q_2}
C_{j_2j_1}\zeta_{j_1}^{(i_1)}\zeta_{j_2}^{(i_2)},
\end{equation}

\vspace{1mm}
\begin{equation}
\label{112}
J^{*}[\psi^{(3)}]_{T,t}=
\hbox{\vtop{\offinterlineskip\halign{
\hfil#\hfil\cr
{\rm l.i.m.}\cr
$\stackrel{}{{}_{q_1,q_2,q_3\to \infty}}$\cr
}} }\sum_{j_1=0}^{q_1}\sum_{j_2=0}^{q_2}\sum_{j_3=0}^{q_3}
C_{j_3j_2j_1}\zeta_{j_1}^{(i_1)}\zeta_{j_2}^{(i_2)}\zeta_{j_3}^{(i_3)},
\end{equation}

\vspace{1mm}
\begin{equation}
\label{feto19000a}
J^{*}[\psi^{(3)}]_{T,t}=
\hbox{\vtop{\offinterlineskip\halign{
\hfil#\hfil\cr
{\rm l.i.m.}\cr
$\stackrel{}{{}_{q\to \infty}}$\cr
}} }
\sum\limits_{j_1, j_2, j_3=0}^{q}
C_{j_3 j_2 j_1}\zeta_{j_1}^{(i_1)}\zeta_{j_2}^{(i_2)}\zeta_{j_3}^{(i_3)},
\end{equation}

\vspace{1mm}
\begin{equation}
\label{feto1900otit}
J^{*}[\psi^{(4)}]_{T,t}=
\hbox{\vtop{\offinterlineskip\halign{
\hfil#\hfil\cr
{\rm l.i.m.}\cr
$\stackrel{}{{}_{q\to \infty}}$\cr
}} }
\sum\limits_{j_1,j_2,j_3,j_4=0}^{q}
C_{j_4 j_3 j_2 j_1}\zeta_{j_1}^{(i_1)}\zeta_{j_2}^{(i_2)}
\zeta_{j_3}^{(i_3)}\zeta_{j_4}^{(i_4)},
\end{equation}

\vspace{5mm}
\noindent
where we assume that
$i_1, i_2, i_3=1,\ldots,m$
in {\rm (\ref{111})--(\ref{feto19000a})} and 
$i_1, \ldots, i_4=0, 1,\ldots,m$ in {\rm (\ref{feto1900otit})}.
Additionally, we assume in
{\rm (\ref{112})} and {\rm (\ref{feto1900otit})} that
$\psi_1(\tau),\ldots,\psi_4(\tau)\equiv 1$.
Another notations are the same as in 
Theorems {\rm 1, 2}.}

Recently, a new approach to the expansion and mean-square 
approximation of iterated Stratonovich stochastic integrals has been obtained
\cite{20a} (Sect.~2.10--2.16), 
\cite{25} (Sect.~5--11), \cite{arxiv-5} (Sect.~13--19), \cite{arxiv-4} (Sect.~7--13),
\cite{new-art-1-xxy} (Sect.~4--9), \cite{new-art-1xxys}.
Let us formulate four theorems that were obtained using this approach.

\vspace{2mm}

{\bf Theorem 4}\ \cite{20a}, \cite{25}, \cite{arxiv-5},  \cite{arxiv-4}, \cite{new-art-1-xxy}.\
{\it Suppose 
that $\{\phi_j(x)\}_{j=0}^{\infty}$ is a complete orthonormal system of 
Legendre polynomials or trigonometric functions in the space $L_2([t, T]).$
Furthermore, let $\psi_1(\tau), \psi_2(\tau),$ $\psi_3(\tau)$ are continuously dif\-ferentiable 
nonrandom functions on $[t, T].$ 
Then, for the 
iterated Stra\-to\-no\-vich stochastic integral of third multiplicity

$$
J^{*}[\psi^{(3)}]_{T,t}={\int\limits_t^{*}}^T\psi_3(t_3)
{\int\limits_t^{*}}^{t_3}\psi_2(t_2)
{\int\limits_t^{*}}^{t_2}\psi_1(t_1)
d{\bf w}_{t_1}^{(i_1)}
d{\bf w}_{t_2}^{(i_2)}d{\bf w}_{t_3}^{(i_3)}\ \ \ (i_1,i_2,i_3=0,1,\ldots,m)
$$

\vspace{4mm}
\noindent
the following 
relations

\vspace{-1mm}
\begin{equation}
\label{fin1}
J^{*}[\psi^{(3)}]_{T,t}
=\hbox{\vtop{\offinterlineskip\halign{
\hfil#\hfil\cr
{\rm l.i.m.}\cr
$\stackrel{}{{}_{p\to \infty}}$\cr
}} }
\sum\limits_{j_1, j_2, j_3=0}^{p}
C_{j_3 j_2 j_1}\zeta_{j_1}^{(i_1)}\zeta_{j_2}^{(i_2)}\zeta_{j_3}^{(i_3)},
\end{equation}

\vspace{3mm}
\begin{equation}
\label{fin2}
{\sf M}\left\{\left(
J^{*}[\psi^{(3)}]_{T,t}-
\sum\limits_{j_1, j_2, j_3=0}^{p}
C_{j_3 j_2 j_1}\zeta_{j_1}^{(i_1)}\zeta_{j_2}^{(i_2)}\zeta_{j_3}^{(i_3)}\right)^2\right\}
\le \frac{C}{p}
\end{equation}

\vspace{5mm}
\noindent
are fulfilled, where $i_1, i_2, i_3=0,1,\ldots,m$ in {\rm (\ref{fin1})} and 
$i_1, i_2, i_3=1,\ldots,m$ in {\rm (\ref{fin2})},
constant $C$ is independent of $p,$

$$
C_{j_3 j_2 j_1}=\int\limits_t^T\psi_3(t_3)\phi_{j_3}(t_3)
\int\limits_t^{t_3}\psi_2(t_2)\phi_{j_2}(t_2)
\int\limits_t^{t_2}\psi_1(t_1)\phi_{j_1}(t_1)dt_1dt_2dt_3
$$

\vspace{4mm}
\noindent
and
$$
\zeta_{j}^{(i)}=
\int\limits_t^T \phi_{j}(\tau) d{\bf f}_{\tau}^{(i)}
$$ 

\vspace{2mm}
\noindent
are independent standard Gaussian random variables for various 
$i$ or $j$ {\rm (}in the case when $i\ne 0${\rm );} 
another notations are the same as in Theorems~{\rm 1, 2}.}

\vspace{2mm}

{\bf Theorem 5}\ \cite{20a}, \cite{25}, \cite{arxiv-5},  \cite{arxiv-4},
\cite{new-art-1-xxy}.\ {\it Let
$\{\phi_j(x)\}_{j=0}^{\infty}$ be a complete orthonormal system of 
Legendre polynomials or trigonometric functions in the space $L_2([t, T]).$
Furthermore, let $\psi_1(\tau), \ldots,$ $\psi_4(\tau)$ be continuously dif\-ferentiable 
nonrandom functions on $[t, T].$ 
Then, for the 
iterated Stra\-to\-no\-vich stochastic integral of fourth multiplicity

\begin{equation}
\label{fin0}
J^{*}[\psi^{(4)}]_{T,t}={\int\limits_t^{*}}^T\psi_4(t_4)
{\int\limits_t^{*}}^{t_4}\psi_3(t_3)
{\int\limits_t^{*}}^{t_3}\psi_2(t_2)
{\int\limits_t^{*}}^{t_2}\psi_1(t_1)
d{\bf w}_{t_1}^{(i_1)}
d{\bf w}_{t_2}^{(i_2)}d{\bf w}_{t_3}^{(i_3)}d{\bf w}_{t_4}^{(i_4)}
\end{equation}

\vspace{4mm}
\noindent
the following 
relations

\begin{equation}
\label{fin3}
J^{*}[\psi^{(4)}]_{T,t}
=\hbox{\vtop{\offinterlineskip\halign{
\hfil#\hfil\cr
{\rm l.i.m.}\cr
$\stackrel{}{{}_{p\to \infty}}$\cr
}} }
\sum\limits_{j_1, j_2, j_3,j_4=0}^{p}
C_{j_4j_3 j_2 j_1}\zeta_{j_1}^{(i_1)}\zeta_{j_2}^{(i_2)}\zeta_{j_3}^{(i_3)}\zeta_{j_4}^{(i_4)},
\end{equation}

\vspace{3mm}

\begin{equation}
\label{fin4}
{\sf M}\left\{\left(
J^{*}[\psi^{(4)}]_{T,t}-
\sum\limits_{j_1, j_2, j_3, j_4=0}^{p}
C_{j_4 j_3 j_2 j_1}\zeta_{j_1}^{(i_1)}\zeta_{j_2}^{(i_2)}\zeta_{j_3}^{(i_3)}
\zeta_{j_4}^{(i_4)}
\right)^2\right\}
\le \frac{C}{p^{1-\varepsilon}}
\end{equation}

\vspace{5mm}
\noindent
are fulfilled, where $i_1, \ldots , i_4=0,1,\ldots,m$ in {\rm (\ref{fin0}),} {\rm (\ref{fin3})} 
and $i_1, \ldots, i_4=1,\ldots,m$ in {\rm (\ref{fin4}),}
constant $C$ does not depend on $p,$
$\varepsilon$ is an arbitrary
small positive real number 
for the case of complete orthonormal system of 
Legendre polynomials in the space $L_2([t, T])$
and $\varepsilon=0$ for the case of
complete orthonormal system of 
trigonometric functions in the space $L_2([t, T]),$

$$
C_{j_4 j_3 j_2 j_1}=
$$

$$
=
\int\limits_t^T\psi_4(t_4)\phi_{j_4}(t_4)
\int\limits_t^{t_4}\psi_3(t_3)\phi_{j_3}(t_3)
\int\limits_t^{t_3}\psi_2(t_2)\phi_{j_2}(t_2)
\int\limits_t^{t_2}\psi_1(t_1)\phi_{j_1}(t_1)dt_1dt_2dt_3dt_4;
$$

\vspace{4mm}
\noindent
another notations are the same as in Theorem~{\rm 4}.}

\vspace{2mm}

{\bf Theorem 6}\ \cite{20a}, \cite{25}, \cite{arxiv-5},  \cite{arxiv-4},
\cite{new-art-1-xxy}.\
{\it Assume 
that $\{\phi_j(x)\}_{j=0}^{\infty}$ is a complete orthonormal system of 
Legendre polynomials or trigonometric functions in the space $L_2([t, T])$
and $\psi_1(\tau), \ldots,$ $\psi_5(\tau)$ are continuously dif\-ferentiable 
nonrandom functions on $[t, T].$ 
Then, for the 
iterated Stra\-to\-no\-vich stochastic integral of fifth multiplicity

\begin{equation}
\label{fin7}
J^{*}[\psi^{(5)}]_{T,t}={\int\limits_t^{*}}^T\psi_5(t_5)
\ldots
{\int\limits_t^{*}}^{t_2}\psi_1(t_1)
d{\bf w}_{t_1}^{(i_1)}
\ldots d{\bf w}_{t_5}^{(i_5)}
\end{equation}

\vspace{4mm}
\noindent
the following 
relations

\begin{equation}
\label{fin8}
J^{*}[\psi^{(5)}]_{T,t}
=\hbox{\vtop{\offinterlineskip\halign{
\hfil#\hfil\cr
{\rm l.i.m.}\cr
$\stackrel{}{{}_{p\to \infty}}$\cr
}} }
\sum\limits_{j_1,\ldots,j_5=0}^{p}
C_{j_5 \ldots j_1}\zeta_{j_1}^{(i_1)}\ldots \zeta_{j_5}^{(i_5)},
\end{equation}

\vspace{3mm}

\begin{equation}
\label{fin9}
{\sf M}\left\{\left(
J^{*}[\psi^{(5)}]_{T,t}-
\sum\limits_{j_1, \ldots, j_5=0}^{p}
C_{j_5 \ldots j_1}\zeta_{j_1}^{(i_1)}\ldots
\zeta_{j_5}^{(i_5)}
\right)^2\right\}
\le \frac{C}{p^{1-\varepsilon}}
\end{equation}

\vspace{5mm}
\noindent
are fulfilled, where $i_1, \ldots , i_5=0,1,\ldots,m$ in {\rm (\ref{fin7}),} {\rm (\ref{fin8})} 
and $i_1, \ldots, i_5=1,\ldots,m$ in {\rm (\ref{fin9}),}
constant $C$ is independent of $p,$
$\varepsilon$ is an arbitrary
small positive real number 
for the case of complete orthonormal system of 
Legendre polynomials in the space $L_2([t, T])$
and $\varepsilon=0$ for the case of
complete orthonormal system of 
trigonometric functions in the space $L_2([t, T]),$

$$
C_{j_5 \ldots j_1}=
\int\limits_t^T\psi_5(t_5)\phi_{j_5}(t_5)\ldots
\int\limits_t^{t_2}\psi_1(t_1)\phi_{j_1}(t_1)dt_1\ldots dt_5;
$$

\vspace{3mm}
\noindent
another notations are the same as in Theorems~{\rm 4, 5}.}

\vspace{2mm}

{\bf Theorem 7}\ \cite{20a}, \cite{25}, \cite{arxiv-5},  \cite{arxiv-4},
\cite{new-art-1xxys}.\
{\it Suppose that 
$\{\phi_j(x)\}_{j=0}^{\infty}$ is a complete orthonormal system of 
Legendre polynomials or trigonometric functions in the space $L_2([t, T]).$
Then, for the 
iterated Stratonovich stochastic integral of sixth multiplicity

\vspace{1mm}
\begin{equation}
\label{after10001qu1}
J_{T,t}^{*(i_1\ldots i_6)}={\int\limits_t^{*}}^T
\ldots
{\int\limits_t^{*}}^{t_2}
d{\bf w}_{t_1}^{(i_1)}
\ldots d{\bf w}_{t_6}^{(i_6)}
\end{equation}

\vspace{5mm}
\noindent
the following 
expansion

$$
J_{T,t}^{*(i_1\ldots i_6)}
=\hbox{\vtop{\offinterlineskip\halign{
\hfil#\hfil\cr
{\rm l.i.m.}\cr
$\stackrel{}{{}_{p\to \infty}}$\cr
}} }
\sum\limits_{j_1, \ldots, j_6=0}^{p}
C_{j_6 \ldots j_1}\zeta_{j_1}^{(i_1)}\ldots
\zeta_{j_6}^{(i_6)}
$$

\vspace{5mm}
\noindent
that converges in the mean-square sense is valid, where
$i_1, \ldots, i_6=0, 1,\ldots,m,$

\vspace{1mm}
$$
C_{j_6 \ldots j_1}=
\int\limits_t^T\phi_{j_6}(t_6)\ldots
\int\limits_t^{t_2}\phi_{j_1}(t_1)dt_1\ldots dt_6;
$$

\vspace{4mm}
\noindent
another notations are the same as in Theorems~{\rm 4--6}.}

\vspace{5mm}

\section{Legendre Polynomial-Based Approximation of the Iterated
Ito and Stratonovich Stochastic Integrals Used in the 
Applications}

\vspace{5mm}

We notice that the collection of iterated Ito stochastic 
integrals used in the numerical
methods (\ref{4.18}), (\ref{4.35}) is given by

\vspace{-1mm}
\begin{equation}
\label{uuu1}
I_{(0)T,t}^{(i_1)},\ \ \ I_{(1)T,t}^{(i_1)},\ \ \ I_{(00)T,t}^{(i_1i_2)},\ \ \
I_{(000)T,t}^{(i_1i_2i_3)},\ \ \ I_{(01)T,t}^{(i_1i_2)},\ \ \
I_{(10)T,t}^{(i_1i_2)},\ \ \ I_{(0000)T,t}^{(i_1i_2i_3i_4)},
\end{equation}

\vspace{3mm}
\noindent
where $i_1,\ldots, i_4=1,\ldots,m$.

The functions $K(t_1,\ldots,t_k)$ like (\ref{ppp}) for the collection
(\ref{uuu1}) are given, respectively, by

\vspace{2mm}
$$
K_0(t_1)\equiv 1,\ \ \ 
K_1(t_1)=t-t_1,\ \ \ K_{00}(t_1,t_2)={\bf 1}_{\{t_1<t_2\}},
$$

$$
K_{000}(t_1,t_2,t_3)={\bf 1}_{\{t_1<t_2<t_3\}},\ \ \
K_{01}(t_1,t_2)=(t-t_2){\bf 1}_{\{t_1<t_2\}},
$$

$$
K_{10}(t_1,t_2)=(t-t_1){\bf 1}_{\{t_1<t_2\}},\ \ \
K_{0000}(t_1,\ldots,t_4)={\bf 1}_{\{t_1<t_2<t_3<t_4\}},
$$

\vspace{5mm}
\noindent
where $t_1,\ldots,t_4\in [t, T]$ and
${\bf 1}_A$ is the indicator of the set $A$.                        

For a finite-degree polynomial, the simplest (having a finite number 
of terms) expansion into
Fourier series by the complete orthonormal system of functions 
in the space 
$L_2([t, T])$ is the Fourier--Legendre series expansion. 
The polynomial functions are 
included in the functions
$K_1(t_1)$, $K_{01}(t_1,t_2)$, $K_{10}(t_1,t_2)$ as their components. 
Therefore, it is logical to expect that the
simplest expansions of these functions 
into multiple Fourier series are their Fourier--Legendre
expansions. 

The following example illustrates rather well the noticed feature.

Consider the approximation $I_{(1)T,t}^{(i_1)q}$
of the stochastic integral $I_{(1)T,t}^{(i_1)}$ 
based on the expansion of
the Brownian bridge process into the trigonometric Fourier series with 
random coefficients \cite{Mi2}

\begin{equation}
\label{42}
I_{(1)T,t}^{(i_1)q}=-\frac{{(T-t)}^{3/2}}{2}
\Biggl(\zeta_0^{(i_1)}-\frac{\sqrt{2}}{\pi}\Biggl(\sum_{r=1}^{q}
\frac{1}{r}
\zeta_{2r-1}^{(i_1)}+\sqrt{\alpha_q}\xi_q^{(i_1)}\Biggr)
\Biggr),
\end{equation}

\vspace{4mm}
\noindent
where
$$
\xi_q^{(i_1)}=\frac{1}{\sqrt{\alpha_q}}\sum_{r=q+1}^{\infty}
\frac{1}{r}\zeta_{2r-1}^{(i_1)},\ \ \
\alpha_q=\frac{\pi^2}{6}-\sum_{r=1}^q\frac{1}{r^2},
$$

\vspace{4mm}
\noindent
where
$\zeta_0^{(i_1)},$ 
$\zeta_{2r-1}^{(i_1)},$ $\xi_q^{(i_1)};$ $r=1,\ldots,q;$
$i_1=1,\ldots,m$ are independent standard Gaussian
random variables.

On the other hand, it is possible to obtain the following equality 

\begin{equation}
\label{4002}
I_{(1)T,t}^{(i_1)}=-\frac{(T-t)^{3/2}}{2}\left(\zeta_0^{(i_1)}+
\frac{1}{\sqrt{3}}\zeta_1^{(i_1)}\right),
\end{equation}

\vspace{4mm}
\noindent
which is valid 
with probability 1 and
based on the expansion of the 
function $t-t_1$ into the Fourier--Legendre series
at the interval $[t,T]$
(this expansion has just two terms). 

The above example 
demonstrates the advantage of the Legendre
polynomials over the trigonometric functions in the context 
of the issue under consideration. More detailed comparison 
can be found in \cite{20a}-\cite{20aaa}, \cite{29}, \cite{30}.

We notice that, as was established in \cite{7}, \cite{8}-\cite{16}, 
\cite{19}-\cite{20aaa}, in the Fourier 
method (Theorem 1) it is also possible 
to use the Haar and Rademacher--Walsh functions (also see Theorem 2). However, 
in \cite{7}, \cite{8}-\cite{16}, \cite{19}-\cite{20aaa}
it was shown that the expansions
of the iterated Ito stochastic integrals (\ref{ito}) of multiplicities 
1 and 2 
obtained with the use of Theorem 1 and systems of Haar and 
Rademacher--Walsh functions are
overcomplicated as compared with their analogues obtained on the 
basis of the Legendre polynomials. In this connection, practical 
application of such expansions is
hindered.

Consider approximations of the remaining stochastic integrals from the 
family (\ref{uuu1}) obtained using Theorems
1, 2 and complete orthonormal system of Legendre 
polynomials in the
space $L_2([t, T]).$ First, we consider approximations of 
stochastic integrals of multiplicities 1
and 2 

\vspace{1mm}
\begin{equation}
\label{4001}
I_{(0)T,t}^{(i_1)}=\sqrt{T-t}\zeta_0^{(i_1)},
\end{equation}

\vspace{2mm}
$$
I_{(00)T,t}^{(i_1 i_2)q}=
I_{(00)T,t}^{*(i_1 i_2)q}-
\frac{1}{2}{\bf 1}_{\{i_1=i_2\}}(T-t),
$$

\vspace{3mm}
\begin{equation}
\label{4004}
I_{(00)T,t}^{*(i_1 i_2)q}=
\frac{T-t}{2}\left(\zeta_0^{(i_1)}\zeta_0^{(i_2)}+\sum_{i=1}^{q}
\frac{1}{\sqrt{4i^2-1}}\left(
\zeta_{i-1}^{(i_1)}\zeta_{i}^{(i_2)}-
\zeta_i^{(i_1)}\zeta_{i-1}^{(i_2)}\right)\right),
\end{equation}

\vspace{3mm}
\begin{equation}
\label{4004a}
I_{(10)T,t}^{(i_1 i_2)q}=
I_{(10)T,t}^{*(i_1 i_2)q}+
\frac{1}{4}{\bf 1}_{\{i_1=i_2\}}(T-t)^2,\ \ \
I_{(01)T,t}^{(i_1 i_2)q}=
I_{(01)T,t}^{*(i_1 i_2)q}+
\frac{1}{4}{\bf 1}_{\{i_1=i_2\}}(T-t)^2,
\end{equation}

\vspace{5mm}
$$
I_{(01)T,t}^{*(i_1 i_2)q}=-\frac{T-t}{2}I_{(00)T,t}^{*(i_1 i_2)q}
-\frac{(T-t)^2}{4}\Biggl(
\frac{1}{\sqrt{3}}\zeta_0^{(i_1)}\zeta_1^{(i_2)}+\Biggr.
$$

\vspace{2mm}
\begin{equation}
\label{4005}
+\Biggl.\sum_{i=0}^{q}\Biggl(
\frac{(i+2)\zeta_i^{(i_1)}\zeta_{i+2}^{(i_2)}
-(i+1)\zeta_{i+2}^{(i_1)}\zeta_{i}^{(i_2)}}
{\sqrt{(2i+1)(2i+5)}(2i+3)}-
\frac{\zeta_i^{(i_1)}\zeta_{i}^{(i_2)}}{(2i-1)(2i+3)}\Biggr)\Biggr),
\end{equation}

\vspace{6mm}

$$
I_{(10)T,t}^{*(i_1 i_2)q}=
-\frac{T-t}{2}I_{(00)T,t}^{*(i_1 i_2)q}
-\frac{(T-t)^2}{4}\Biggl(
\frac{1}{\sqrt{3}}\zeta_0^{(i_2)}\zeta_1^{(i_1)}+\Biggr.
$$

\vspace{2mm}
\begin{equation}
\label{4006}
+\Biggl.\sum_{i=0}^{q}\Biggl(
\frac{(i+1)\zeta_{i+2}^{(i_2)}\zeta_{i}^{(i_1)}
-(i+2)\zeta_{i}^{(i_2)}\zeta_{i+2}^{(i_1)}}
{\sqrt{(2i+1)(2i+5)}(2i+3)}+
\frac{\zeta_i^{(i_1)}\zeta_{i}^{(i_2)}}{(2i-1)(2i+3)}\Biggr)\Biggr),
\end{equation}

\vspace{7mm}
\noindent
where here and below

$$
I_{(l_1\ldots l_k)s,t}^{*(i_1\ldots i_k)q}\ \ \ \hbox{and}\ \ \
I_{(l_1\ldots l_k)s,t}^{(i_1\ldots i_k)q}
$$

\vspace{5mm}
\noindent 
are the approximations of the iterated Stratonovich
and Ito stochastic integrals like 

\begin{equation}
\label{502.1aye}
I_{(l_1\ldots l_k)s,t}^{*(i_1\ldots i_k)}=
{\int\limits_t^{*}}^s (t-\tau _
{k}) ^ {l_ {k}}\ldots 
{\int\limits_t^{*}}^{\tau_{2}} (t-\tau _ {1}) ^ {l_ {1}} d{\bf f} ^ {(i_ {1})} _ {\tau_ {1}} \ldots 
d {\bf f} _ {\tau_ {k}} ^ {(i_ {k})}
\end{equation}

\vspace{4mm}
\noindent
and, correspondingly, like (\ref{ll1}); $\zeta_{j}^{(i)}$
are independent standard Gaussian
random variables for various $i$ or $j$;
$j=0, 1,\ldots,p+2;$\ 
$i=1,\ldots,m.$

Calculate the mean-square errors of approximations
(\ref{4004})--(\ref{4006}). A precise formula for pairwise
different $i_1,\ldots,i_k=1,\ldots,m$ was established in 
\cite{7}, \cite{20}-\cite{20aaa}, \cite{26}

\vspace{1mm}
\begin{equation}
\label{102}
{\sf M}\left\{\left(J[\psi^{(k)}]_{T,t}-J[\psi^{(k)}]_{T,t}^q\right)^2\right\}=
\int\limits_{[t,T]^k} K^2(t_1,\ldots,t_k)dt_1\ldots dt_k
-\sum_{j_1,\ldots,j_k=0}^q C_{j_k\ldots j_1}^2,
\end{equation}

\vspace{4mm}
\noindent
where in virtue of the Parseval equality (\ref{sos2z})
the right-hand side of (\ref{102})
tends to zero for $q\to\infty;$ $J[\psi^{(k)}]_{T,t}$ has
the form (\ref{ito}), and $J[\psi^{(k)}]_{T,t}^q$
is the approximation of $J[\psi^{(k)}]_{T,t}$ 
defined as the prelimit expression
in (\ref{leto6000}) for $p_1=\ldots=p_k=q$ (also see the prelimit 
expressions in (\ref{a1})--(\ref{a5})); the sense of the
rest notations is the same as in Theorems 1, 2.

The following formula \cite{7}, \cite{20}-\cite{20aaa}, \cite{26} takes place

\vspace{1mm}
$$
{\sf M}\left\{\left(J[\psi^{(2)}]_{T,t}-J[\psi^{(2)}]_{T,t}^q\right)^2\right\}
=
$$

\vspace{1mm}
\begin{equation}
\label{at1}
=\int\limits_{[t,T]^2}K^2(t_1,t_2)dt_1dt_2-\sum_{j_1,j_2=0}^q
C_{j_2j_1}^2-\sum_{j_1,j_2=0}^q C_{j_1j_2}C_{j_2j_1}\ \ \ (i_1=i_2),
\end{equation}

\vspace{5mm}
\noindent
where notations are the same as in (\ref{102}).

The value ${\sf M}\left\{\left(J[\psi^{(k)}]_{T,t}-J[\psi^{(k)}]_{T,t}^q\right)^2\right\}$
can be calculated exactly.

\vspace{2mm}

{\bf Theorem 8} \cite{20a} (Sect.~1.12), \cite{26} (Sect.~6).
{\it Suppose that $\{\phi_j(x)\}_{j=0}^{\infty}$ 
is an arbitrary complete orthonormal system  
of functions in the space $L_2([t,T])$ and
$\psi_1(\tau),\ldots,\psi_k(\tau)\in L_2([t, T]),$  $i_1,\ldots, i_k=1,\ldots,m$.
Then

$$
{\sf M}\left\{\left(J[\psi^{(k)}]_{T,t}-J[\psi^{(k)}]_{T,t}^q\right)^2\right\}=
\int\limits_{[t,T]^k} K^2(t_1,\ldots,t_k)dt_1\ldots dt_k- 
$$

\begin{equation}
\label{tttr11}
-
\sum_{j_1,\ldots, j_k=0}^{q}
C_{j_k\ldots j_1}
{\sf M}\left\{J[\psi^{(k)}]_{T,t}
\sum\limits_{(j_1,\ldots,j_k)}
\int\limits_t^T \phi_{j_k}(t_k)
\ldots
\int\limits_t^{t_{2}}\phi_{j_{1}}(t_{1})
d{\bf f}_{t_1}^{(i_1)}\ldots
d{\bf f}_{t_k}^{(i_k)}\right\},
\end{equation}

\vspace{5mm}
\noindent
where $i_1,\ldots,i_k = 1,\ldots,m;$
the expression 

\vspace{-1mm}
$$
\sum\limits_{(j_1,\ldots,j_k)}
$$ 

\vspace{3mm}
\noindent
means the sum with respect to all
possible permutations 
$(j_1,\ldots,j_k)$. At the same time if 
$j_r$ swapped with $j_q$ in the permutation $(j_1,\ldots,j_k),$
then $i_r$ swapped with $i_q$ in the permutation
$(i_1,\ldots,i_k);$
another notations are the same as in Theorems {\rm 1, 2.}
}

\vspace{2mm}

Using (\ref{102}) and (\ref{at1}), we get

\vspace{1mm}
\begin{equation}
\label{4007}
{\sf M}\left\{\left(I_{(00)T,t}^{(i_1 i_2)}-
I_{(00)T,t}^{(i_1 i_2)q}
\right)^2\right\}
=\frac{(T-t)^2}{2}\Biggl(\frac{1}{2}-\sum_{i=1}^q
\frac{1}{4i^2-1}\Biggr)\ \ \ (i_1\ne i_2),
\end{equation}

\vspace{6mm}
$$
{\sf M}\left\{\left(I_{(10)T,t}^{(i_1 i_2)}-
I_{(10)T,t}^{(i_1 i_2)q}
\right)^2\right\}=
{\sf M}\left\{\left(I_{(01)T,t}^{(i_1 i_2)}-
I_{(01)T,t}^{(i_1 i_2)q}\right)^2\right\}=\frac{(T-t)^4}{16}\times
$$

\vspace{2mm}
\begin{equation}
\label{2007ura}
\times\left(\frac{5}{9}-
2\sum_{i=2}^q\frac{1}{4i^2-1}-
\sum_{i=1}^q
\frac{1}{(2i-1)^2(2i+3)^2}
-\sum_{i=0}^q\frac{(i+2)^2+(i+1)^2}{(2i+1)(2i+5)(2i+3)^2}
\right)
\end{equation}

\vspace{6mm}
\noindent
for $i_1\ne i_2$ and

\vspace{2mm}
$$
{\sf M}\left\{\left(I_{(10)T,t}^{(i_1 i_1)}-
I_{(10)T,t}^{(i_1 i_1)q}
\right)^2\right\}=
{\sf M}\left\{\left(I_{(01)T,t}^{(i_1 i_1)}-
I_{(01)T,t}^{(i_1 i_1)q}\right)^2\right\}=
$$

\vspace{2mm}
\begin{equation}
\label{2007ura1}
=\frac{(T-t)^4}{16}\left(\frac{1}{9}-
\sum_{i=0}^{q}
\frac{1}{(2i+1)(2i+5)(2i+3)^2}-
2\sum_{i=1}^{q}
\frac{1}{(2i-1)^2(2i+3)^2}\right).
\end{equation}

\vspace{6mm}

Let us consider the numerical modeling of the iterated 
Ito stochastic integral of multiplicity 3
$I_{(000)T,t}^{(i_1i_2i_3)}.$ Using 
Theorems 1, 2 for the case $k=3$ (see (\ref{a3})), we obtain

\vspace{1mm}
$$
I_{(000)T,t}^{(i_1i_2i_3)q}
=\sum_{j_1,j_2,j_3=0}^{q}
C_{j_3j_2j_1}\Biggl(
\zeta_{j_1}^{(i_1)}\zeta_{j_2}^{(i_2)}\zeta_{j_3}^{(i_3)}
-{\bf 1}_{\{i_1=i_2\}}
{\bf 1}_{\{j_1=j_2\}}
\zeta_{j_3}^{(i_3)}-
\Biggr.
$$

\vspace{2mm}
\begin{equation}
\label{sad001}
\Biggl.
-{\bf 1}_{\{i_2=i_3\}}
{\bf 1}_{\{j_2=j_3\}}
\zeta_{j_1}^{(i_1)}-
{\bf 1}_{\{i_1=i_3\}}
{\bf 1}_{\{j_1=j_3\}}
\zeta_{j_2}^{(i_2)}\Biggr),
\end{equation}

\vspace{5mm}
\noindent
where $i_1, i_2, i_3=1,\ldots,m$ and

$$
C_{j_3j_2j_1}=\int\limits_t^T\phi_{j_3}(z)
\int\limits_t^{z} \phi_{j_2}(y)
\int\limits_t^{y}
\phi_{j_1}(x) dx dy dz=
$$

\vspace{2mm}
\begin{equation}
\label{sad00222}
=\frac{\sqrt{(2j_1+1)(2j_2+1)(2j_3+1)}}{8}(T-t)^{3/2}\bar
C_{j_3j_2j_1},
\end{equation}

\vspace{2mm}
\begin{equation}
\label{sad0023}
\bar C_{j_3j_2j_1}=\int\limits_{-1}^{1}P_{j_3}(z)
\int\limits_{-1}^{z}P_{j_2}(y)
\int\limits_{-1}^{y}
P_{j_1}(x)dx dy dz,
\end{equation}

\vspace{5mm}
\noindent
where $P_i(x)$ $(i=0, 1, 2,\ldots)$ is the Legendre polynomial.

For the case $i_1=i_2=i_3$, one can use the well known 
equality which follows from the Ito
formula and is valid with probability 1 \cite{KlPl2}

\begin{equation}
\label{sad003}
I_{(000)T,t}^{(i_1 i_1 i_1)}
=\frac{1}{6}(T-t)^{3/2}\left(
\left(\zeta_0^{(i_1)}\right)^3-3
\zeta_0^{(i_1)}\right).
\end{equation}

\vspace{4mm}

The procedure of numerical modeling of the iterated Ito
stochastic integral 
$I_{(000)T,t}^{(i_1i_2i_3)}$ may follow (\ref{sad001})--(\ref{sad003}).
The Fourier--Legendre coefficients $\bar C_{j_3j_2j_1}$ of the form 
(\ref{sad0023}) being precisely calculable 
for the given number $q$ by PYTHON, DERIVE or MAPLE.
The mean-square error of approximation is
checked by (\ref{102}) for $k=3$
as well as by the formulas established in 
\cite{20}-\cite{20aaa}, \cite{26}

\vspace{1mm}
$$
{\sf M}\left\{\left(J[\psi^{(3)}]_{T,t}-
J[\psi^{(3)}]_{T,t}^q\right)^2\right\}
=\int\limits_{[t,T]^3}K^2(t_1,t_2,t_3)dt_1dt_2dt_3-
$$

\vspace{1mm}
\begin{equation}
\label{2000}
-\sum_{j_3,j_2,j_1=0}^q C_{j_3j_2j_1}^2-
\sum_{j_3,j_2,j_1=0}^q C_{j_3j_1j_2}C_{j_3j_2j_1}\ \ \ (i_1=i_2\ne i_3),
\end{equation}

\vspace{7mm}
$$
{\sf M}\left\{\left(J[\psi^{(3)}]_{T,t}-
J[\psi^{(3)}]_{T,t}^q\right)^2\right\}
=\int\limits_{[t,T]^3}K^2(t_1,t_2,t_3)dt_1dt_2dt_3-
$$

\vspace{1mm}
\begin{equation}
\label{2001}
-
\sum_{j_3,j_2,j_1=0}^q C_{j_3j_2j_1}^2-
\sum_{j_3,j_2,j_1=0}^q C_{j_2j_3j_1}C_{j_3j_2j_1}\ \ \ (i_1\ne i_2=i_3),
\end{equation}

\vspace{7mm}
$$
{\sf M}\left\{\left(J[\psi^{(3)}]_{T,t}-
J[\psi^{(3)}]_{T,t}^q\right)^2\right\}
=\int\limits_{[t,T]^3}K^2(t_1,t_2,t_3)dt_1dt_2dt_3-
$$

\vspace{1mm}
\begin{equation}
\label{2002}
-\sum_{j_3,j_2,j_1=0}^q C_{j_3j_2j_1}^2-
\sum_{j_3,j_2,j_1=0}^q C_{j_3j_2j_1}C_{j_1j_2j_3}\ \ \ (i_1=i_3\ne i_2).
\end{equation}

\vspace{7mm}
\noindent

The following estimate \cite{20}-\cite{20aaa}, \cite{26}
can also be applied for the case $k=3$

\vspace{2mm}
$$
{\sf M}\left\{\left(J[\psi^{(k)}]_{T,t}-
J[\psi^{(k)}]_{T,t}^q\right)^2\right\}
\le 
$$

\begin{equation}
\label{tu1}
\le k!\left(
\int\limits_{[t,T]^k}
K^2(t_1,\ldots,t_k)
dt_1\ldots dt_k -\sum_{j_1,\ldots,j_k=0}^{q}C^2_{j_k\ldots j_1}\right)
\end{equation}

\vspace{5mm}
\noindent
where $i_1,\ldots,i_k=1,\ldots,m$ and $0<T-t<\infty$ or
$i_1,\ldots,i_k=0, 1,\ldots,m$ and $0<T-t<1.$

In particular, for the pairwise different 
$i_1,i_2,i_3=1,\ldots,m$ and $q=6$ we get from (\ref{102})

\vspace{-2mm}
\begin{equation}
\label{rrrr2}
{\sf M}\left\{\left(
I_{(000)T,t}^{(i_1i_2 i_3)}-
I_{(000)T,t}^{(i_1i_2 i_3)6}\right)^2\right\}\approx
0.01956(T-t)^3.
\end{equation}

\vspace{3mm}

Taking into consideration that $T-t$ is
the integration step of numerical
methods for the Ito SDE (\ref{1.5.2}) and $T-t$
is a sufficiently small number, we get that already
for $q=6$ the mean-square error of approximation of the 
stochastic integral $I_{{000}_{T,t}}^{(i_1i_2 i_3)}$ is sufficiently
small as well (see (\ref{rrrr2})).

Consider now the iterated Ito stochastic integral 
$I_{(0000)T,t}^{(i_1 i_2 i_3 i_4)}$
of multiplicity 4. Using Theorems 1, 2, we get
the representation

$$
I_{(0000)T,t}^{(i_1 i_2 i_3 i_4)q}
=\sum_{j_1,j_2,j_3,j_4=0}^{q}
C_{j_4 j_3 j_2 j_1}\Biggl(
\zeta_{j_1}^{(i_1)}\zeta_{j_2}^{(i_2)}\zeta_{j_3}^{(i_3)}\zeta_{j_4}^{(i_4)}
-\Biggr.
$$
$$
-
{\bf 1}_{\{i_1=i_2\}}
{\bf 1}_{\{j_1=j_2\}}
\zeta_{j_3}^{(i_3)}
\zeta_{j_4}^{(i_4)}
-
{\bf 1}_{\{i_1=i_3\}}
{\bf 1}_{\{j_1=j_3\}}
\zeta_{j_2}^{(i_2)}
\zeta_{j_4}^{(i_4)}-
$$
$$
-
{\bf 1}_{\{i_1=i_4\}}
{\bf 1}_{\{j_1=j_4\}}
\zeta_{j_2}^{(i_2)}
\zeta_{j_3}^{(i_3)}
-
{\bf 1}_{\{i_2=i_3\}}
{\bf 1}_{\{j_2=j_3\}}
\zeta_{j_1}^{(i_1)}
\zeta_{j_4}^{(i_4)}-
$$
$$
-
{\bf 1}_{\{i_2=i_4\}}
{\bf 1}_{\{j_2=j_4\}}
\zeta_{j_1}^{(i_1)}
\zeta_{j_3}^{(i_3)}
-
{\bf 1}_{\{i_3=i_4\}}
{\bf 1}_{\{j_3=j_4\}}
\zeta_{j_1}^{(i_1)}
\zeta_{j_2}^{(i_2)}+
$$
$$
+
{\bf 1}_{\{i_1=i_2\}}
{\bf 1}_{\{j_1=j_2\}}
{\bf 1}_{\{i_3=i_4\}}
{\bf 1}_{\{j_3=j_4\}}
+
{\bf 1}_{\{i_1=i_3\}}
{\bf 1}_{\{j_1=j_3\}}
{\bf 1}_{\{i_2=i_4\}}
{\bf 1}_{\{j_2=j_4\}}+
$$
\begin{equation}
\label{399}
+\Biggl.
{\bf 1}_{\{i_1=i_4\}}
{\bf 1}_{\{j_1=j_4\}}
{\bf 1}_{\{i_2=i_3\}}
{\bf 1}_{\{j_2=j_3\}}\Biggr),
\end{equation}

\vspace{6mm}
\noindent
where $i_1, i_2, i_3, i_4=1,\ldots,m$ and

\vspace{2mm}
$$
C_{j_4j_3j_2j_1}=
\int\limits_{t}^{T}\phi_{j_4}(u)
\int\limits_{t}^{u}\phi_{j_3}(z)
\int\limits_{t}^{z} \phi_{j_2}(y)
\int\limits_{t}^{y}
\phi_{j_1}(x) dx dy dz du=
$$

\vspace{3mm}
$$
=\frac{\sqrt{(2j_1+1)(2j_2+1)(2j_3+1)(2j_4+1)}}{16}\Delta^{2}\bar
C_{j_4j_3j_2j_1},
$$

\vspace{4mm}
$$
\bar C_{j_4j_3j_2j_1}=\int\limits_{-1}^{1}P_{j_4}(u)
\int\limits_{-1}^{u}P_{j_3}(z)
\int\limits_{-1}^{z}P_{j_2}(y)
\int\limits_{-1}^{y}
P_{j_1}(x)dx dy dz du,
$$

\vspace{5mm}
\noindent
where $P_i(x)$ $(i=0, 1, 2,\ldots)$ is the Legendre polynomial.

For precise calculation of the Fourier--Legendre coefficients 
$C_{j_4 j_3 j_2 j_1}$
we can use the previous recommendations
and check the mean-square error of approximation of the 
iterated Ito stochastic integral
$I_{(0000)T,t}^{(i_1i_2i_3 i_4)}$, 
for example, using the estimate (\ref{tu1}) for $k=4$.

In particular, for pairwise different 
$i_1,\ldots,i_4=1,\ldots,m$ we get from (\ref{102})
with regard for
smallness of $T-t$ already for $q=2$ a sufficiently good 
accuracy of the mean-square approximation

\begin{equation}
\label{rrrr1}
{\sf M}\left\{\left(
I_{(0000)T,t}^{(i_1i_2i_3 i_4)}-
I_{(0000)T,t}^{(i_1i_2i_3 i_4)2}\right)^2\right\}\approx
0.0236084(T-t)^4.
\end{equation}

\vspace{4mm}

We notice that at deriving (\ref{rrrr2}) and (\ref{rrrr1}) 
the coefficients $\bar C_{j_3j_2j_1}$ and $\bar C_{j_4j_3j_2j_1}$ 
were precisely
calculated using the DERIVE package.

Note that the formulas (\ref{111})--(\ref{feto1900otit}) 
are simpler than (\ref{a2})--(\ref{a4}). However,
calculation of the mean-square approximation error 
for the iterated Stratonovich stochastic integrals 
(\ref{str}) turned out more complex than for the iterated Ito
stochastic integrals 
(\ref{ito}) \cite{20a}-\cite{20aaa}, \cite{28}, \cite{99999}.

\vspace{5mm}

\section{Algorithms of Numerical Modeling With the Orders
$1.5$ and $2.0$ of Strong Convergence}

\vspace{5mm}

We formulate in algorithmic form the above formulas and 
recommendations for the numerical
method of the order $1.5$ of strong convergence. We assume 
that the necessary Fourier--Legendre coefficients
$\bar C_{j_3j_2j_1},$
$\bar C_{j_4j_3j_2j_1}$
are already calculated. In particular, 
several tables of the precisely calculated Fourier--Legendre coefficients
$\bar C_{j_3j_2j_1},$ $\bar C_{j_4j_3j_2j_1}$ 
were presented in 
\cite{7}, \cite{20}-\cite{20aaa}. These coefficients
were calculated by DERIVE.
It should be noted that in \cite{Kuz-Kuz}, \cite{Mikh-1} the database with
270,000 precisely calculated Fourier--Legendre coefficients
is presented. In \cite{Kuz-Kuz}, \cite{Mikh-1} 
we used the PYTHON programming language.

\vspace{5mm}

\centerline{\bf Algorithm. 1.}

\vspace{5mm}

{\bf Step 1.}\ Given are the initial parameters of the problem 
such as the interval of integration $[0, T]$,
step of integration $\Delta$ (for example, constant $\Delta=T/N,$ $N\ge 1$,
although a variable step of integration
is admissible), initial condition ${\bf y}_0$, and constant $C$ involved 
in the condition (\ref{4.3}).

\vspace{4mm}

{\bf Step 2.}\ Assume that $p=0.$

\vspace{4mm}

{\bf Step 3.}\ 
Selection of the minimal natural numbers $q$ and $q_1$ ($q\ll q_1$) 
ensuring the necessary accuracy
of approximation of the stochastic integrals 

\vspace{-1mm}
$$
I_{(00)\tau_{p+1},\tau_p}^{(i_1 i_2)},\ \ \
I_{(000)\tau_{p+1},\tau_k}^{(i_1 i_2 i_3)}\ \ \ (\tau_p=p\Delta)
$$

\vspace{3mm}
\noindent
and satisfying the conditions

\begin{equation}
\label{sad0091x}
{\sf M}\left\{\left(I_{(00)\tau_{p+1},\tau_p}^{(i_1 i_2)}-
{I}_{(00)\tau_{p+1},\tau_p}^{(i_1 i_2)q_1}
\right)^2\right\}
=\frac{\Delta^2}{2}\Biggl(\frac{1}{2}-\sum_{i=1}^{q_1}
\frac{1}{4i^2-1}\Biggr)\le C\Delta^4,
\end{equation}

\vspace{2mm}
\begin{equation}
\label{sad0091}
{\sf M}\left\{\left(
I_{(000)\tau_{p+1},\tau_p}^{(i_1i_2 i_3)}-
I_{(000)\tau_{p+1},\tau_p}^{(i_1i_2 i_3)q}\right)^2\right\}\le
6\left(\frac{\Delta^{3}}{6}-\sum_{j_3,j_2,j_1=0}^{q}
C_{j_3j_2j_1}^2\right)\le C\Delta^4.
\end{equation}

\vspace{5mm}

{\bf Remark 1.}\
{\it If it is required to check the mean-square approximation error
of the iterated Ito stochastic 
integral $I_{(000)\tau_{p+1},\tau_p}^{(i_1i_2 i_3)}$
using the
precise formulas {\rm (\ref{102}), (\ref{2000})--(\ref{2002})}, 
rather than the 
estimate {\rm (\ref{sad0091})} {\rm (}see {\rm (\ref{tu1}))}, 
then instead of
the condition {\rm (\ref{sad0091})} one has to take the 
following conditions

\vspace{1mm}
$$
E^{(i_1i_2i_3)}_{p,q,\Delta}=
\frac{\Delta^{3}}{6}-\sum_{j_3,j_2,j_1=0}^{q}
C_{j_3j_2j_1}^2\le C\Delta^4\ \ \ (i_1\ne i_2,\ i_1\ne i_3,\ i_2\ne i_3),
$$

\vspace{3mm}
$$
E^{(i_1i_2i_3)}_{p,q,\Delta}=
\frac{\Delta^{3}}{6}-\sum_{j_3,j_2,j_1=0}^{q}
C_{j_3j_2j_1}^2
-\sum_{j_3,j_2,j_1=0}^{q}
C_{j_2j_3j_1}C_{j_3j_2j_1}\le C\Delta^4\ \ \ (i_1\ne i_2=i_3),
$$

\vspace{3mm}
$$
E^{(i_1i_2i_3)}_{p,q,\Delta}=
\frac{\Delta^{3}}{6}-\sum_{j_3,j_2,j_1=0}^{q}
C_{j_3j_2j_1}^2
-\sum_{j_3,j_2,j_1=0}^{q}
C_{j_3j_2j_1}C_{j_1j_2j_3}\le C\Delta^4\ \ \ (i_1=i_3\ne i_2),
$$

\vspace{3mm}
$$
E^{(i_1i_2i_3)}_{p,q,\Delta}=
\frac{\Delta^{3}}{6}-\sum_{j_3,j_2,j_1=0}^{q}
C_{j_3j_2j_1}^2-\sum_{j_3,j_2,j_1=0}^{q}
C_{j_3j_1j_2}C_{j_3j_2j_1}\le C\Delta^4\ \ \ (i_1=i_2\ne i_3),
$$

\vspace{6mm}
\noindent
where

\vspace{-1mm}
$$
{\sf M}\left\{\left(
I_{(000)\tau_{p+1},\tau_p}^{(i_1i_2 i_3)}-
I_{(000)\tau_{p+1},\tau_p}^{(i_1i_2 i_3)q}\right)^2\right\}
\stackrel{{\rm def}}{=}E^{(i_1i_2i_3)}_{p,q,\Delta}.
$$

\vspace{5mm}

{\bf Step 4.}\ Modeling of the sequence of independent 
standard Gaussian random variables
$\zeta_l^{(i)}$ $(l=0, 1,\ldots,q_1;\ i=1,\ldots,m).$

\vspace{4mm}

{\bf Step 5.}\ Modeling of the iterated Ito stochastic integrals

$$
I_{(0)\tau_{p+1},\tau_p}^{(i_1)},\ \ \
I_{(1)\tau_{p+1},\tau_p}^{(i_1)},\ \ \
I_{(00)\tau_{p+1},\tau_p}^{(i_1 i_2)},\ \ \ 
I_{(000)\tau_{p+1},\tau_p}^{(i_1 i_2 i_3)}
$$

\vspace{4mm}
\noindent
using the formulas

\vspace{-2mm}
$$
I_{(0)\tau_{k+1},\tau_k}^{(i_1)}=\sqrt{T-t}\zeta_0^{(i_1)},
$$

\vspace{3mm}
$$
I_{(1)\tau_{p+1},\tau_p}^{(i_1)}=-\frac{(T-t)^{3/2}}{2}\left(\zeta_0^{(i_1)}+
\frac{1}{\sqrt{3}}\zeta_1^{(i_1)}\right),
$$

\vspace{5mm}
$$
I_{(00)\tau_{p+1},\tau_p}^{(i_1 i_2)q_1}=
\frac{T-t}{2}\left(\zeta_0^{(i_1)}\zeta_0^{(i_2)}+\sum_{i=1}^{q_1}
\frac{1}{\sqrt{4i^2-1}}\left(
\zeta_{i-1}^{(i_1)}\zeta_{i}^{(i_2)}-
\zeta_i^{(i_1)}\zeta_{i-1}^{(i_2)}\right)-{\bf 1}_{\{i_1=i_2\}}\right),
$$

\vspace{7mm}

$$
I_{(000)\tau_{p+1},\tau_p}^{(i_1i_2i_3)q}
=\sum_{j_1,j_2,j_3=0}^{q}
C_{j_3j_2j_1}\Biggl(
\zeta_{j_1}^{(i_1)}\zeta_{j_2}^{(i_2)}\zeta_{j_3}^{(i_3)}
-{\bf 1}_{\{i_1=i_2\}}
{\bf 1}_{\{j_1=j_2\}}
\zeta_{j_3}^{(i_3)}-
\Biggr.
$$

\vspace{3mm}
$$
\Biggl.
-{\bf 1}_{\{i_2=i_3\}}
{\bf 1}_{\{j_2=j_3\}}
\zeta_{j_1}^{(i_1)}-
{\bf 1}_{\{i_1=i_3\}}
{\bf 1}_{\{j_1=j_3\}}
\zeta_{j_2}^{(i_2)}\Biggr),
$$

\vspace{4mm}
\noindent
where $i_1, i_2, i_3 =1,\ldots,m.$

\vspace{2mm}

{\bf Remark 2.}\ {\it In the case of $i_1=i_2=i_3$, it is advisable 
to model the stochastic integral $I_{(000)\tau_{p+1},\tau_p}^{(i_1i_2i_3)}$
using the formula {\rm (\ref{sad003}),} 
where one has to assume that $T-t=\Delta.$}

\vspace{4mm}

{\bf Step 6.}\ Calculate ${\bf y}_{p+1}$ from {\rm (\ref{4.18})}.

\vspace{4mm}

{\bf Step 7.}\ If $p<N-1$, then assume that $p=p+1$ and 
go to Step {\rm 4;} otherwise, go to Step {\rm 8}.

\vspace{4mm}

{\bf Step 8.}\ End.}

\vspace{4mm}

We briefly note how to modify the algorithm to enable numerical 
modeling with the order $2.0$ of
strong convergence. 

At Step 3 one has to take the 
following three iterated Ito stochastic integrals 

$$
I_{(10)\tau_{p+1},\tau_p}^{(i_1i_2)},\ \ \
I_{(01)\tau_{p+1},\tau_p}^{(i_1i_2)},\ \ \
I_{(0000)\tau_{p+1},\tau_p}^{(i_1i_2 i_3 i_4)},
$$

\vspace{3mm}
\noindent
whose approximations obey (\ref{4004a})--(\ref{4006}), (\ref{399})
and add to the considered stochastic 
integrals. Moreover, we replace $C\Delta^4$ by $C\Delta^5$
in (\ref{sad0091x}), (\ref{sad0091}).
At that, one can use
the estimate (\ref{tu1}) for $k=4$ and the formulas (\ref{2007ura}),
(\ref{2007ura1}) to 
check the accuracy of modeling of the
aforementioned integrals. As the result, we get the following conditions

\vspace{2mm}
$$
{\sf M}\left\{\left(I_{(10)\tau_{p+1},\tau_p}^{(i_1 i_2)}-
I_{(10)\tau_{p+1},\tau_p}^{(i_1 i_2)q_2}
\right)^2\right\}=
{\sf M}\left\{\left(I_{(01)\tau_{p+1},\tau_p}^{(i_1 i_2)}-
I_{(01)\tau_{p+1},\tau_p}^{(i_1 i_2)q_2}\right)^2\right\}=
\frac{\Delta^4}{16}\times
$$

\vspace{3mm}
$$
\times
\left(\frac{5}{9}-
2\sum_{i=2}^{q_2}\frac{1}{4i^2-1}-
\sum_{i=1}^{q_2}
\frac{1}{(2i-1)^2(2i+3)^2}-
\sum_{i=0}^{q_2}\frac{(i+2)^2+(i+1)^2}{(2i+1)(2i+5)(2i+3)^2}
\right)\le C\Delta^5
$$

\vspace{6mm}
\noindent
for $i_1\ne i_2$ and

\vspace{2mm}
$$
{\sf M}\left\{\left(I_{(10)\tau_{p+1},\tau_p}^{(i_1 i_1)}-
I_{(10)\tau_{p+1},\tau_p}^{(i_1 i_1)q_3}
\right)^2\right\}=
{\sf M}\left\{\left(I_{(01)\tau_{p+1},\tau_p}^{(i_1 i_1)}-
I_{(01)\tau_{p+1},\tau_p}^{(i_1 i_1)q_3}\right)^2\right\}=
$$

\vspace{4mm}
$$
=\frac{\Delta^4}{16}\left(\frac{1}{9}-
\sum_{i=0}^{q_3}
\frac{1}{(2i+1)(2i+5)(2i+3)^2}-
2\sum_{i=1}^{q_3}
\frac{1}{(2i-1)^2(2i+3)^2}\right)\le C\Delta^5
$$

\vspace{6mm}
\noindent
for $i_1=i_2;$

\vspace{2mm}
\begin{equation}
\label{jj1}
{\sf M}\left\{\left(
I_{(0000)\tau_{p+1},\tau_p}^{(i_1i_2 i_3 i_4)}-
I_{(0000)\tau_{p+1},\tau_p}^{(i_1i_2 i_3 i_4)q_4}\right)^2\right\}\le
24\left(\frac{\Delta^{4}}{24}-\sum_{j_1,j_2,j_3,j_4=0}^{q_4}
C_{j_4j_3j_2j_1}^2\right)\le C\Delta^5,
\end{equation}

\vspace{5mm}
\noindent
where $i_1, i_2, i_3, i_4=1,\ldots,m;$\ $q_2, q_3, q_4<q<q_1.$

Carry out Step 5 with allowance of the stochastic integrals 

$$
I_{(10)\tau_{p+1},\tau_p}^{(i_1 i_2)},\ 
I_{(01)\tau_{p+1},\tau_p}^{(i_1 i_2)},\
I_{(0000)\tau_{p+1},\tau_p}^{(i_1 i_2 i_3 i_4)},
$$

\vspace{4mm}
\noindent
and calculate ${\bf y}_{p+1}$ at Step 6 according to (\ref{4.35}).
 
It should be noted that instead of the estimate (\ref{jj1})
we can use the precise relations for the value

$$
{\sf M}\left\{\left(
I_{(0000)\tau_{p+1},\tau_p}^{(i_1i_2 i_3 i_4)}-
I_{(0000)\tau_{p+1},\tau_p}^{(i_1i_2 i_3 i_4)q_4}\right)^2\right\},
$$

\vspace{4mm}
\noindent
which were obtained in \cite{20}-\cite{20aaa}, \cite{26}
for all possible combinations
of $i_1,i_2,i_3,i_4=1,\ldots,m$.
Note that the optimization of the mentioned procedure is
considered in \cite{Mikh-2}.

\vspace{5mm}

\section{Conclusions}

\vspace{5mm}

The present paper provided efficient procedures for the mean-square 
approximation of iterated Ito and
Stratonovich stochastic integrals of multiplicities 
1 to 4 based on multiple Fourier--Legendre series. 
These results can be used for implementation of the numerical 
methods with the orders $1.0$, $1.5$, and $2.0$ 
of strong convergence for Ito stochastic 
differential equations with multidimensional 
non-commutative noise.
The results of the article can be applied 
for numerical solution
of the problems of optimal stochastic control 
and signal filtering in random noise in different
formulations. The development of the approaches from this 
work can be found in 
\cite{7}, \cite{8}-\cite{31}, \cite{Kuz-Kuz}-\cite{art-6qqqqq}, \cite{Mikh-2},
\cite{new-new-1}.


\vspace{12mm}

\begin{thebibliography}{199}

\vspace{8mm}


\bibitem{1}
Gihman I.I., Skorohod A.V. Stochastic Differential Equations and 
its Applications.
Naukova Dumka, Kiev, 1982, 612 pp.


\bibitem{KlPl2}
Kloeden P.E., Platen E. Numerical Solution of Stochastic
Differential Equations. Springer, Berlin, 1995, 632 pp.


\bibitem{KPS}
Kloeden P.E., Platen E., Schurz H. Numerical Solution of SDE 
Through Computer 
Experiments. Springer, Berlin, 1994, 292 pp.



\bibitem{Arato}
Arato M. Linear Stochastic Systems With Constant Coefficients. 
A Statistical Approach. Springer, Berlin, Heidelberg, N.Y., 
1982, 289 pp. 


\bibitem{Shir1}
Shiriaev A.N. Foundations of Financial Mathematics. Vol. 2,
Fazis, Moscow, 1998, 544 pp. 


\bibitem{Lip}
Liptser R.Sh., Shiriaev A.N. Statistics of Stochastic Processes:
Nonlinear Filtering and Related Problems. Nauka, Moscow, 1974,
696 pp. 


\bibitem{Mi2}
Milstein G.N. Numerical Integration of Stochastic Differential 
Equations. Ural University Press, Sverdlovsk, 1988, 225 pp.  


\bibitem{Mi3}
Milstein G.N., Tretyakov M.V. 
Stochastic Numerics for Mathematical Physics. 
Springer, Berlin, 2004, 616 pp.



\bibitem{PW1}
Platen E., Wagner W. On a Taylor formula for a class of Ito
processes. Probab. Math. Statist. 3 (1982), 37-51.


\bibitem{KlPl1}
Kloeden P.E., Platen E. The Stratonovich and Ito-Taylor 
expansions. Math. Nachr. 151 (1991), 33-50.


\bibitem{kk5}
Kulchitskiy O.Yu., Kuznetsov D.F. The unified Taylor-Ito expansion.
Journal of Mathematical Sciences (N. Y.) 99, 2 (2000), 1130-1140.	
DOI: https://doi.org/10.1007/BF02673635


\bibitem{kuz33}
Kuznetsov, D.F. New representations of the Taylor-Stratonovich expansion.
Journal of Mathematical Sciences (N. Y.). 118, 6 (2003), 5586-5596.\
DOI: http://doi.org/10.1023/A:1026138522239



\bibitem{7}
Kuznetsov D.F. Numerical Integration of Stochastic Differential Equations. 2.
[In Russian]. Polytechnical University Publishing House, 
Saint-Petersburg, 2006, 764 pp.
DOI: http://doi.org/10.18720/SPBPU/2/s17-227\
Available at:\ http://www.sde-kuznetsov.spb.ru/06.pdf\
(ISBN 5-7422-1191-0)



\bibitem{3}
Kuznetsov D.F. A method of expansion and approximation of repeated 
stochastic Stratonovich integrals based on multiple Fourier series 
on full orthonormal systems. [In Russian].
Electronic Journal "Differential Equations and Control Processes"
ISSN 1817-2172 (online),
1 (1997), 18-77.
Available at:\\
http://diffjournal.spbu.ru/EN/numbers/1997.1/article.1.2.html



\bibitem{4}
Kuznetsov D.F. Problems of the Numerical Analysis of Ito Stochastic 
Differential Equations. 
[In Russian].
Electronic Journal "Differential Equations and Control Processes"
ISSN 1817-2172 (online),
1 (1998), 66-367.
Available at:\
http://diffjournal.spbu.ru/EN/numbers/1998.1/article.1.3.html\
Hard Cover Edition: 1998, SPbGTU Publishing House, 
204 pp. (ISBN 5-7422-0045-5)


\bibitem{KPW}
Kloeden P.E., Platen E., Wright I.W. The approximation of multiple 
stochastic integrals. Stochastic Analysis and Applications.
10, 4 (1992), 431-441. 


\bibitem{Zapad-9}
Platen, E., Bruti-Liberati N. Numerical Solution of Stochastic 
Differential Equations
with Jumps in Finance. Springer, Berlin-Heidelberg, 2010, 868 pp.



\bibitem{Al}
Allen E. Approximation of triple stochastic integrals through region 
subdivision. Communications in Applied Analysis
(Special Tribute Issue to Professor V. Lakshmikantham). 17 (2013), 355-366. 


\bibitem{8}
Kuznetsov D.F. Stochastic Differential Equations: Theory and Practice 
of Numerical Solution. With MatLab Programs, 1st Edition. [In Russian]. 
Polytechnical University Publishing House, Saint-Petersburg, 2007, 778 pp.
DOI: http://doi.org/10.18720/SPBPU/2/s17-228\
Available at:\ http://www.sde-kuznetsov.spb.ru/07b.pdf\
(ISBN 5-7422-1394-8)




\bibitem{9}
Kuznetsov D.F. Stochastic Differential Equations: Theory and Practice 
of Numerical Solution. With MatLab Programs, 2nd Edition. [In Russian]. 
Polytechnical University Publishing House, Saint-Petersburg, 2007, XXXII+770 pp.
DOI: http://doi.org/10.18720/SPBPU/2/s17-229\
Available at:\\
http://www.sde-kuznetsov.spb.ru/07a.pdf\
(ISBN 5-7422-1439-1)




\bibitem{10}
Kuznetsov D.F. Stochastic Differential Equations: Theory and Practice 
of Numerical Solution. With MatLab Programs, 3rd Edition. [In Russian]. 
Polytechnical University Publishing House, Saint-Petersburg, 2009, XXXIV+768 pp.
DOI: http://doi.org/10.18720/SPBPU/2/s17-230\
Available at:\\
http://www.sde-kuznetsov.spb.ru/09.pdf\
(ISBN 978-5-7422-2132-6)





\bibitem{11}
Kuznetsov D.F. Stochastic Differential Equations: Theory and Practice 
of Numerical Solution. With MatLab Programs. 4th Edition. [In Russian].
Polytechnical University Publishing House, Saint-Petersburg, 2010, 
XXX+786 pp. DOI: http://doi.org/10.18720/SPBPU/2/s17-231\
Available at:\
http://www.sde-kuznetsov.spb.ru/10.pdf\
(ISBN 978-5-7422-2448-8)



\bibitem{12}
Kuznetsov D.F. Multiple Stochastic Ito and Stratonovich Integrals 
and Multiple Fourier Series.
[In Russian].
Electronic Journal "Differential Equations and Control Processes"
ISSN 1817-2172 (online),
3 (2010), A.1-A.257. DOI: http://doi.org/10.18720/SPBPU/2/z17-7\
Available at:\\
http://diffjournal.spbu.ru/EN/numbers/2010.3/article.2.1.html





\bibitem{13}
Kuznetsov D.F. Strong Approximation of Multiple Ito and 
Stratonovich Stochastic Integrals: Multiple Fourier Series Approach.
1st Edition. [In English]. 
Polytechnical University Publishing House, Saint-Petersburg, 
2011, 250 pp. DOI: http://doi.org/10.18720/SPBPU/2/s17-232\
Available at:\\
http://www.sde-kuznetsov.spb.ru/11b.pdf\
(ISBN 978-5-7422-2988-9)



\bibitem{14}
Kuznetsov D.F. Strong Approximation of Multiple Ito and 
Stratonovich Stochastic Integrals: Multiple Fourier Series Approach.
2nd Edition. [In English]. 
Polytechnical University Publishing House, Saint-Petersburg, 
2011, 284 pp. DOI: http://doi.org/10.18720/SPBPU/2/s17-233\
Available at:\\
http://www.sde-kuznetsov.spb.ru/11a.pdf\
(ISBN 978-5-7422-3162-2)



\bibitem{15}
Kuznetsov D.F. Multiple Ito and Stratonovich Stochastic 
Integrals: Approximations, Properties, Formulas. [In English].
Polytechnical University Publishing House, Saint-Petersburg,
2013, 382 pp.\\
DOI: http://doi.org/10.18720/SPBPU/2/s17-234\\
Available at:\
http://www.sde-kuznetsov.spb.ru/13.pdf\
(ISBN 978-5-7422-3973-4)




\bibitem{16}
Kuznetsov D.F. Multiple Ito and Stratonovich Stochastic Integrals: 
Fourier-Legendre and Trigonometric Expansions, Approximations, Formulas.
[In English].
Electronic Journal "Differential Equations and Control Processes"
ISSN 1817-2172 (online),
1 (2017), A.1--A.385.\
DOI: http://doi.org/10.18720/SPBPU/2/z17-3\\ 
Available at:\
http://diffjournal.spbu.ru/EN/numbers/2017.1/article.2.1.html 


\bibitem{17}
Kuznetsov D.F. Development and application of the Fourier 
method for the numerical solution of Ito stochastic differential 
equations. [In English]. Computational Mathematics and 
Mathematical Physics, 58, 7 (2018), 1058-1070.
DOI: http://doi.org/10.1134/S0965542518070096



\bibitem{18}
Kuznetsov D.F. On numerical modeling of the multidimensional 
dynamic systems under random perturbations with the 1.5 and 2.0 
orders of strong convergence [In English]. Automation and Remote Control,
79, 7 (2018), 1240-1254.
DOI: http://doi.org/10.1134/S0005117918070056


\bibitem{19}
Kuznetsov D.F. Stochastic Differential Equations: Theory and Practice of 
Numerical Solution. With Programs on MATLAB, 5th Edition. [In Russian].
Electronic Journal "Differential Equations and Control Processes"
ISSN 1817-2172 (online), 2 (2017), A.1-A.1000.
DOI: http://doi.org/10.18720/SPBPU/2/z17-4\
Available at:\\
http://diffjournal.spbu.ru/EN/numbers/2017.2/article.2.1.html




\bibitem{20}
Kuznetsov D.F. Stochastic Differential Equations: Theory and Practice of 
Numerical Solution. With MATLAB Programs, 6th Edition. [In Russian].
Electronic Journal "Differential Equations and Control Processes"
ISSN 1817-2172 (online), 4 (2018), A.1-A.1073.\
Available at:\\
http://diffjournal.spbu.ru/EN/numbers/2018.4/article.2.1.html


\bibitem{20a}
Kuznetsov D.F.
Strong Approximation of Iterated Ito and Stratonovich Stochastic 
Integrals Based on Generalized Multiple Fourier Series. 
Application to Numerical Solution of Ito SDEs and Semilinear SPDEs.
[In English].
arXiv:2003.14184 [math.PR]. 2022, 923 pp.


\bibitem{20aa}
Kuznetsov D.F.
Strong Approximation of Iterated It\^{o} and Stratonovich Stochastic 
Integrals Based on Generalized Multiple Fourier Series. 
Application to Numerical Solution of It\^{o} SDEs and Semilinear SPDEs.
[In English].
Electronic Journal "Differential Equations and Control Processes"
ISSN 1817-2172 (online),
4 (2020), A.1-A.606.\ Available at:\
http://diffjournal.spbu.ru/EN/numbers/2020.4/article.1.8.html



\bibitem{20aaa}
Kuznetsov D.F.
Mean-Square Approximation of Iterated It\^{o} and Stratonovich Stochastic 
Integrals Based on Generalized Multiple Fourier Series. 
Application to Numerical Integration of It\^{o} SDEs and Semilinear SPDEs.
[In English].
Electronic Journal "Differential Equations and Control Processes"
ISSN 1817-2172 (online),
4 (2021), A.1-A.788.\ Available at:\
http://diffjournal.spbu.ru/EN/numbers/2021.4/article.1.9.html




\bibitem{22}
Kuznetsov D.F. Mean-square approximation of iterated Ito 
and Stratonovich stochastic 
integrals of multiplicities 1 to 6 
from the Taylor-Ito and Taylor-Stratonovich expansions
using Legendre polynomials. 
[In English]. arXiv:1801.00231 [math.PR]. 2019, 106 pp. 




\bibitem{23}
Kuznetsov D.F.
Expansion of iterated Stratonovich stochastic integrals 
of arbitrary multiplicity
based on generalized iterated Fourier series converging pointwise. 
[In English].
arXiv:1801.00784 [math.PR]. 2018, 77 pp. 



\bibitem{25}
Kuznetsov D.F.
Expansion of iterated Stratonovich stochastic integrals of fifth and sixth
multiplicity based on generalized multiple Fourier series.
[In English]. arXiv:1802.00643 [math.PR]. 2022, 129 pp. 



\bibitem{26}
Kuznetsov D.F.
Exact calculation of the mean-square error in the method
of approximation of iterated Ito stochastic 
integrals based on generalized multiple Fourier series. [In English].
arXiv:1801.01079 [math.PR]. 2018, 68 pp. 


\bibitem{26a}
Kuznetsov D.F. 
Expansion of iterated Ito stochastic integrals of arbitrary multiplicity
based on generalized multiple Fourier series converging in the mean. 
[In English]. arXiv:1712.09746 [math.PR]. 2022, 111 pp.


\bibitem{26b}
Kuznetsov D.F. Expansion of iterated stochastic integrals with respect
to martingale Poisson measures and with respect to martingales based on 
generalized multiple Fourier series. [In English]. arXiv:1801.06501 [math.PR].
2018, 40 pp.




\bibitem{27}
Kuznetsov D.F. Expansion of iterated Stratonovich stochastic integrals
based on generalized multiple Fourier series. 
[In English]. Ufa Mathematical Journal, 
11, 4 (2019), 49-77. 
DOI: http://doi.org/10.13108/2019-11-4-49\\
Available at:\
http://matem.anrb.ru/en/article?art\_id=604.




\bibitem{28}
Kuznetsov D.F. On numerical modeling of the multidimentional dynamic 
systems under random perturbations with the 2.5 order of strong 
convergence. [In English]. Automation and Remote Control, 
80, 5 (2019), 867-881. DOI: http://doi.org/10.1134/S0005117919050060




\bibitem{29}
Kuznetsov D.F. Comparative analysis of the efficiency of application 
of Legendre polynomials and trigonometric functions to the numerical 
integration of It\^{o} stochastic differential equations.
[In English]. Computational Mathematics and Mathematical Physics, 
59, 8 (2019),  1236-1250.\\
DOI: http://doi.org/10.1134/S0965542519080116




\bibitem{30}
Kuznetsov D.F. Comparative analysis of the efficiency of application 
of Legendre polynomials and trigonometric functions to the numerical 
integration of Ito stochastic differential equations.
[In English]. arXiv:1901.02345 [math.GM], 2019, 40 pp.


\bibitem{30a}
Kuznetsov D.F. Expansion of multiple Stratonovich stochastic integrals 
of second multiplicity based on double Fourier-Legendre series 
summarized by Prinsheim method [In Russian]. 
Electronic Journal "Differential Equations and Control Processes"
ISSN 1817-2172 (online), 1 (2018), 1-34.\
Available at:\\
http://diffjournal.spbu.ru/EN/numbers/2018.1/article.1.1.html
 

\bibitem{31a}
Kuznetsov D.F. Application of the method of approximation 
of iterated stochastic Ito integrals based on generalized 
multiple Fourier series to the high-order strong numerical methods 
for non-commutative semilinear stochastic partial 
differential equations [In English]. Electronic Journal "Differential 
Equations and Control Processes" 
ISSN 1817-2172 (online), 3 (2019), 18-62.\
Available at:\\
http://diffjournal.spbu.ru/EN/numbers/2019.3/article.1.2.html



\bibitem{arxiv-5}
Kuznetsov D.F. 
Expansions of iterated Stratonovich stochastic integrals
based on generalized multiple Fourier series:
multiplicities 1 to 6 and beyond. [in English].
arXiv:1712.09516 [math.PR]. 2022, 204 pp. 





\bibitem{arxiv-8}
Kuznetsov D.F.  
Expansion of iterated Stratonovich stochastic integrals of 
multiplicity 2 based on double Fourier-Legendre series 
summarized by Pringsheim method. [in English].
arXiv:1801.01962 [math.PR]. 2018, 49 pp.




\bibitem{arxiv-4}
Kuznetsov D.F.
The hypotheses on expansions of iterated Stratonovich stochastic integrals 
of arbitrary multiplicity and their partial proof. [in English].
arXiv:1801.03195 [math.PR]. 
2022, 138 pp. 


\bibitem{arxiv-44}
Kuznetsov D.F. Numerical simulation of 2.5-set of iterated
Stratonovich stochastic integrals of multiplicities 1 to 5
from the Taylor--Stratovovich expansion. [in English].
arXiv: 1806.10705 
[math.PR]. 2018, 29 pp. 


\bibitem{iii}
Kuznetsov D.F.
Strong numerical methods of orders 2.0, 2.5, and 3.0 for 
Ito stochastic differential equations based on the unified 
stochastic Taylor expansions and multiple Fourier-Legendre series.
[in English].
arXiv:1807.02190 [math.PR], 2018, 44 pp.



\bibitem{99999}
Kuznetsov D.F. Explicit one-step strong numerical methods 
of order 2.0 and 2.5 for Ito stochastic differential 
equations based on the unified Taylor-Ito and Taylor-Stratonovich 
Expansions. [in English].
arXiv: 1802.04844 [math.PR]. 2018, 37 pp.


\bibitem{arxiv-7}
Kuznetsov D.F. Expansion of iterated Stratonovich stochastic integrals
of multiplicity 3 based on generalized 
multiple Fourier series converging in the mean: 
general case of series summation. [in English].
arXiv:1801.01564 [math.PR]. 2018, 65 pp. 




\bibitem{arxiv-12}
Kuznetsov D.F. Development and application of the Fourier method to the 
mean-square approximation 
of iterated Ito and Stratonovich stochastic integrals. [in English].
arXiv:1712.08991 [math.PR]. 2017, 57 pp. 




\bibitem{31}
Kuznetsov D.F. Application of the method of approximation of iterated 
Ito stochastic integrals based on generalized multiple Fourier 
series to the high-order strong numerical methods for non-commutative 
semilinear stochastic partial differential equations. [In English].
arXiv:1905.03724 [math.GM], 2019, 41 pp.




\bibitem{W-Z-1}
Wong E., Zakai M. On the convergence of ordinary integrals to 
stochastic integrals. Ann. Math. Stat.,
5, 36 (1965), 1560-1564.


\bibitem{W-Z-2}
Wong E., Zakai M. On the relation between ordinary and stochastic 
differential equations. Int. J. Eng. Sci., 3 (1965), 213-229.


\bibitem{Watanabe}
Ikeda N., Watanabe S. Stochastic
Differential Equations and Diffusion Processes.
2nd Edition. North-Holland Publishing Company,
Amsterdam, Oxford, New-York, 1989. 555 pp.


\bibitem{5}
Kuznetsov D.F. Mean square approximation of solutions 
of stochastic differential 
equations using Legendres polynomials. [In English]. Journal of 
Automation and 
Information Sciences (Begell House), 2000, 32 (Issue 12), 69-86.
DOI: http://doi.org/10.1615/JAutomatInfScien.v32.i12.80\\
Available at:\ 
http://www.sde-kuznetsov.spb.ru/00a.pdf 



\bibitem{6}
Kuznetsov D.F. New representations of explicit one-step numerical 
methods for jump-diffusion stochastic differential equations. 
[In English]. Computational 
Mathematics and Mathematical Physics, 41, 6 (2001), 874-888.
Available at:\ 
http://www.sde-kuznetsov.spb.ru/01b.pdf




\bibitem{Kuz-Kuz}
Kuznetsov M.D., Kuznetsov D.F.
SDE-MATH: A software package for the implementation of strong high-order 
numerical methods for Ito SDEs with multidimensional non-commutative noise 
based on multiple Fourier--Legendre series. [In English].
Electronic Journal "Differential 
Equations and Control Processes" 
ISSN 1817-2172 (online),
1 (2021), 93-422. Available at:\\
http://diffjournal.spbu.ru/EN/numbers/2021.1/article.1.5.html




\bibitem{Mikh-1}
Kuznetsov M.D., Kuznetsov D.F.
Implementation of strong numerical methods 
of orders 0.5, 1.0, 1.5, 2.0, 2.5, and 3.0 for Ito SDEs with non-commutative 
noise based on the unified Taylor-Ito and Taylor-Stratonovich 
Expansions and multiple Fourier-Legendre series.
[In English].
arXiv:2009.14011 [math.PR], 2020, 343 pp. 


\bibitem{OK1000}
Kuznetsov D.F. The proof of convergence with probability 1 
in the method of expansion 
of iterated Ito stochastic integrals based on generalized multiple 
Fourier series.
[In English]. Electronic Journal "Differential 
Equations and Control Processes"
ISSN 1817-2172 (online),
2 (2020), 89-117. Available at:\\ 
http://diffjournal.spbu.ru/RU/numbers/2020.2/article.1.6.html


\bibitem{OK}
Kuznetsov, D.F. Application of multiple Fourier-Legendre series 
to strong exponential Milstein and 
Wagner-Platen methods for non-commutative semilinear stochastic 
partial differential equations.
Electronic Journal "Differential 
Equations and Control Processes"
ISSN 1817-2172 (online),
3 (2020), 129-162.
Available at:\\
http://diffjournal.spbu.ru/RU/numbers/2020.3/article.1.6.html




\bibitem{art-8eeeeeee}
Kuznetsov D.F. Strong approximation of iterated Ito and Stratonovich 
stochastic integrals. Abstracts of talks given at the 4th
International Conference on Stochastic 
Methods (Divnomorskoe, Russia, June 2-9, 2019), 
Theory of Probability and its Applications, 65, 1 (2020), 141-142.\
DOI: http://doi.org/10.1137/S0040585X97T989878


\bibitem{Kuzh-1}
Kuznetsov D.F. 
Application of multiple Fourier-Legendre series to the implementation 
of strong exponential Milstein and Wagner-Platen methods for 
non-commutative semilinear SPDEs. 
Proceedings of the XIII International Conference on Applied 
Mathematics and Mechanics in the Aerospace Industry (AMMAI-2020). 
MAI, Moscow, 2020, pp. 451-453.\
Available at:\ http://www.sde-kuznetsov.spb.ru/20e.pdf


\bibitem{art-6qqqqq}
Kuznetsov D.F. Explicit one-step mumerical method with the 
strong convergence order of 2.5 for Ito stochastic differential 
equations with a multi-dimensional nonadditive noise 
based on the Taylor-Stratonovich expansion.
[In English]. Computational Mathematics and Mathematical Physics, 
60, 3 (2020), 379-389.\\ 
DOI: http://doi.org/10.1134/S0965542520030100



\bibitem{Rybakov1000}
Rybakov K.A. Orthogonal expansion of multiple It\^{o} stochastic integrals.
Electronic Journal "Differential Equations and Control Processes"
ISSN 1817-2172 (online),
3 (2021), 109-140. Available at:\\
http://diffjournal.spbu.ru/EN/numbers/2021.3/article.1.8.html



\bibitem{Mikh-2}
Kuznetsov M.D., Kuznetsov D.F.
Optimization of the mean-square approximation procedures for 
iterated Ito stochastic integrals of multiplicities 1 to 5 from 
the unified Taylor--Ito expansion based on multiple Fourier--Legendre series.
[In English].
arXiv:2010.13564 [math.PR], 2020, 63 pp. 



\bibitem{new-new-1}
Kuznetsov D.F., Kuznetsov M.D. Mean-square approximation of iterated 
stochastic integrals from strong exponential Milstein and Wagner--Platen 
methods for non-commutative semilinear SPDEs based on multiple 
Fourier--Legendre series. Recent Developments in Stochastic Methods and 
Applications. ICSM-5 2020. 
Springer Proceedings in Mathematics \& Statistics, vol. 371, Eds. 
Shiryaev A.N., Samouylov K.E., Kozyrev D.V.
Springer, Cham, 2021, pp. 17-32.\ DOI: http://doi.org/10.1007/978-3-030-83266-7\_2


\bibitem{new-art-1-xxy}
Kuznetsov D.F. A new approach to the series expansion of iterated 
Stratonovich stochastic integrals of arbitrary multiplicity with respect 
to components of the multidimensional Wiener process. [In English].
Electronic Journal "Differential Equations and Control Processes"
ISSN 1817-2172 (online), 2 (2022), 83-186. 
Available at:\\
http://diffjournal.spbu.ru/EN/numbers/2022.2/article.1.6.html




\bibitem{new-art-1xxys}
Kuznetsov D.F. A new approach to the series expansion of iterated 
Stratonovich stochastic integrals of arbitrary multiplicity with 
respect to components of the multidimensional Wiener process. II.
[In English].
Electronic Journal "Differential Equations and Control Processes"
ISSN 1817-2172 (online), 4 (2022). To appear.
Available at:\
http://diffjournal.spbu.ru/EN/collection.html


\bibitem{new-new-7}
Kuznetsov D.F. The three-step strong numerical methods of the orders of 
accuracy 1.0 and 1.5 for Ito stochastic differential equations. [In English]. 
Journal of Automation and Information Sciences (Begell House), 2002, 
34 (Issue 12), 14 pp.\
DOI: http://doi.org/10.1615/JAutomatInfScien.v34.i12.30\\
Available at:\ http://www.sde-kuznetsov.spb.ru/02a.pdf


\bibitem{new-new-8}
Kuznetsov D.F. Finite-difference strong numerical methods of order 
1.5 and 2.0 for stochastic differential Ito equations with 
nonadditive multidimensional noise. [In English].
Journal of Automation and Information Sciences (Begell House), 2001, 33 (Issue 5-8), 13 pp.
DOI: http://doi.org/10.1615/JAutomatInfScien.v33.i5-8.180\\
Available at:\ http://www.sde-kuznetsov.spb.ru/01c.pdf


\end{thebibliography}
\end{document}